\documentclass[10pt]{amsart}
\usepackage{amscd}
\usepackage{amssymb}

\newtheorem{thm}[equation]{Theorem}
\newtheorem{cor}[equation]{Corollary}
\newtheorem{prop}[equation]{Proposition}
\newtheorem{lem}[equation]{Lemma}

\theoremstyle{definition}
\newtheorem{dfn}[equation]{Definition}
\newtheorem{rem}[equation]{Remark}
\newtheorem{exa}[equation]{Example}

\newtheorem{que}[equation]{Question}
\numberwithin{equation}{section}

\newcommand{\iso}{\stackrel{\simeq}{\rightarrow}}

\newcommand{\surj}{\twoheadrightarrow}
\newcommand{\ar}{\rightarrow}
\newcommand{\opn}{\operatorname}
\newcommand{\cat}[1]{\operatorname{\mathsf{#1}}}

\newcommand{\rmitem}[1]{\item[\text{\textup{(#1)}}]}

\newcommand{\mfrak}[1]{\mathfrak{#1}}
\newcommand{\mcal}[1]{\mathcal{#1}}
\newcommand{\msf}[1]{\mathsf{#1}}

\newcommand{\mrm}[1]{\mathrm{#1}}
\newcommand{\mbb}[1]{\mathbb{#1}}

\newcommand{\tup}[1]{\textup{#1}}
\newcommand{\bsym}[1]{\boldsymbol{#1}}

\title{Rings with Auslander Dualizing Complexes}
\author{Amnon Yekutieli and James J. Zhang }
\date{8 September 1998}
\address{A. Yekutieli: Department of Theoretical Mathematics, The
Weizmann Institute of Science, Rehovot 76100, Israel}
\email{amnon@wisdom.weizmann.ac.il}
\address{J.J. Zhang: Department of Mathematics, Box 354350,
University of Washington, Seattle, Washington 98195, USA}
\email{zhang@math.washington.edu}
\thanks{Both authors were supported by the US-Israel Binational
Science Foundation. The second
author was also supported by the NSF and Sloan Research Fellowship}
\subjclass{Primary: 16D90; Secondary: 16E30, 16L30, 16P60}

\begin{document}

\begin{abstract}
A ring with an Auslander dualizing complex is a
generalization of an Auslander-Gorenstein ring.
We show that many results which hold for Auslander-Gorenstein rings
also hold in the more general setting. On the other hand we give
criteria for existence of Auslander dualizing complexes which
show these occur quite frequently.

The most powerful tool we use is the Local Duality Theorem for
connected graded algebras over a field. Filtrations allow the
transfer of results to non-graded algebras.

We also prove some results of a categorical nature, most notably the
functoriality of rigid dualizing complexes.
\end{abstract}

\maketitle
\tableofcontents

\setcounter{section}{-1}
\section{Introduction}

{\em Dualizing complexes} play an essential role in the
Serre-Grothendieck Duality Theory on schemes (see \cite{RD}). The
duality formalism was generalized to noncommutative rings by the first
author, in order to answer some questions which arose in this context,
mainly regarding local duality for noncommutative graded algebras (see
\cite{Ye1}). A version of Serre duality for noncommutative projective
schemes was established using dualizing complexes (see J{\o}rgensen's
\cite{Jo2} and our \cite{YZ1}). Dualizing complexes, and more
generally derived categories, are powerful tools for proving
abstract properties of noncommutative rings. For examples, consider
the noncommutative graded versions of the Auslander-Buchsbaum
Theorem, The Bass Theorem and the No-Holes Theorem for Bass numbers
(see \cite{Jo1} Theorems 3.2, 4.5, 4.6 and 4.8). Under the
synonym ``cotilting complexes,'' dualizing complexes were studied by
Miyachi \cite{Mi2}. Cotilting bimodules occur often in papers on
representations of finite-dimensional algebras, cf.\ \cite{Ha}.

In this paper we will provide further evidence that the dualizing
complex (Definition \ref{a1.1}) is an effective tool for studying
noncommutative rings. We are especially interested in those dualizing
complexes which satisfy an extra homological condition called the {\em
Auslander property} (Definition \ref{a2.1}). The basic idea here is
that if a statement holds for Auslander-Gorenstein rings in the sense
of \cite{Bj}, then an appropriate version of the statement should hold
for rings with Auslander dualizing complexes. The Gorenstein
condition, i.e.\ the ring itself having finite injective dimension, is
considered to be very restrictive; in contrast, having an Auslander
dualizing complex is considered to be a mild condition.

There are a few ways to show existence of Auslander dualizing
complexes. For example, if $A$ is a connected graded algebra (over a
field $k$) with enough normal elements, then $A$ has an Auslander
dualizing complex. Recall that a connected graded $k$-algebra $A$ has
{\em enough normal elements} if every graded prime factor
$A / \mfrak{p} \neq k$ has a nonzero normal element of positive
degree. This class of rings has been studied recently by many
algebraists, because of developments in quantum groups and
noncommutative algebraic geometry.

In this paper we prove a diverse collection of results, whose common
thread is that their proofs are based on the existence of an Auslander
dualizing complex. Throughout $k$ denotes a fixed base field, and a
$k$-algebra means an associative algebra with $1$.

First we generalize \cite{GL} Theorem 1.6, by dropping the Gorenstein
condition. Note that the hypothesis on $\opn{gr} A$ in next theorem
is easy to check in practice (see Example \ref{a8.26}), and we suspect
it can even be weakened.

\begin{thm} \label{a0.1}
Assume $A$ is a normally separated filtered $k$-algebra such that
$\opn{gr} A$ is a noetherian connected graded $k$-algebra with enough
normal elements. Then $\opn{Spec} A$ is catenary.
\end{thm}

This is proved after Corollary \ref{a8.23}.

We generalize some results in \cite{ASZ}, two of which are:

\begin{thm} \label{a0.3}
Assume $A$ is a noetherian $k$-algebra with an Auslander dualizing
complex. Then there is a step duality between the category
$\cat{Mod}_{\mrm{f}} A$ of finitely generated left $A$-modules and
the category $\cat{Mod}_{\mrm{f}} A^{\circ}$ of finitely generated
right $A$-modules.
\end{thm}

Actually we prove a more general result involving two
algebras -- see Theorem \ref{a2.9}. It follows that if the algebra
$A$ has an Auslander dualizing complex then the left and right
Krull dimensions of $A$ are finite.

\begin{thm} \label{a0.10}
Let $A$ be an Auslander-Gorenstein noetherian $k$-algebra.
Assume $A$ has a filtration such that $\opn{gr} A$ is an
AS-Gorenstein noetherian connected graded $k$-algebra
\tup{(}e.g.\ if $A$ is connected graded\tup{)}.
Then $A$ has an artinian self-injective ring of fractions.
\end{thm}

For more details see Theorem \ref{a8.24}.

The notion of characteristic variety $\opn{Ch} M$ of a module $M$ was
introduced in $\mcal{D}$-module theory (cf.\ \cite{Co} p.\ 98). The
next theorem generalizes a result of Gabber by dropping
the Gorenstein condition. The possibility of making this
generalization was suggested by Van den Bergh.

\begin{thm}[Purity of Characteristic Variety] \label{a0.4}
Assume $A$ is a filtered $k$-algebra such that $\opn{gr} A$ is a
commutative connected graded affine $k$-algebra. If $M$ is a finitely
generated $\opn{GKdim}$-pure left $A$-module, then
the characteristic variety $\opn{Ch} M$ is pure.
\end{thm}

Given an algebra $A$ there might exist non-isomorphic Auslander
dualizing complexes over it (see Example \ref{a2.2}(c)). For various
reasons (like functoriality) it is desirable to find dualizing
complexes which are canonical in some sense. In the graded case, the
{\em balanced dualizing complex} (Definition \ref{a6.3}), introduced
by the first author in \cite{Ye1}, is a natural choice. In the
ungraded case, one should consider the {\em rigid dualizing complex}
(Definition \ref{a4.1}) introduced by Van den Bergh \cite{VdB}. Rigid
dualizing complexes (and balanced dualizing complexes in graded case)
are uniquely determined up to a unique isomorphism.
A balanced dualizing complex is always rigid (see \cite{VdB}
Proposition 8.2(2) and our Corollary \ref{a8.6}).
The next theorem on the functoriality of rigid dualizing complexes
is a combination of Theorem \ref{a4.2} and Corollary \ref{a4.6}.
A homomorphism $A \ar B$ of algebras is called {\em finite} if $B$ is
a finitely generated left and right $A$-module.

\begin{thm}
Let $A$ be a noetherian $k$-algebra. A rigid dualizing complex
$R_{A}$ over $A$ is unique up to a unique isomorphism.
If $A \ar B$ is a finite homomorphism and $R_{A}, R_{B}$ are rigid
dualizing complexes over $A, B$ respectively, then there is at most
one morphism
$\opn{Tr}_{B / A} : R_{B} \ar R_{A}$
compatible with the rigidity.
\end{thm}

The existence of a dualizing complex is not automatic.
A very effective criterion for existence of balanced complexes is
given in \cite{VdB} Theorem 6.3 (which is Theorem \ref{a6.5} here).
Van den Bergh's idea was to first prove the Local Duality Theorem,
and then to show this duality is represented by a balanced
dualizing complex.
The Rees algebra allows us to transfer results on graded algebras to
non-graded algebras (see Theorem \ref{a8.3}).

We prove that Auslander rigid dualizing complexes exist for
a large class of rings. First, the similar submodule condition on
graded $A$-modules (see Definition \ref{a7.10}) enables induction
on $\opn{GKdim}$. Combining Theorems \ref{a7.11} and \ref{a7.12}
we get:

\begin{thm} \label{a0.6}
Assume $A$ is a noetherian connected graded $k$-algebra which has a
balanced dualizing complex $R$ and satisfies the
similar submodule condition \tup{(}e.g.\ $A$ is FBN or
has enough normal elements\tup{)}. Then the balanced dualizing
complex $R$ is graded Auslander.
\end{thm}

Next, Theorem \ref{a8.3} says:

\begin{thm} \label{a0.5}
Suppose $A$ is a filtered $k$-algebra such that
the associated graded algebra $\opn{gr} A$ is noetherian.
If $\opn{gr} A$ has a graded Auslander balanced dualizing complex,
then $A$ has an Auslander rigid dualizing complex.
\end{thm}

By results of Grothendieck, a commutative affine connected graded
$k$-algebra has a graded Auslander balanced dualizing complex (this
also follows from Theorem \ref{a0.6}). Since factor rings of universal
enveloping algebras of finite dimensional Lie algebras are filtered,
and their associated graded algebras are commutative, Theorem
\ref{a0.5} tells us that these algebras have Auslander rigid dualizing
complexes. In the same way we may use Theorems \ref{a0.5} and
\ref{a0.6} to show that many quantum algebras and their factor
algebras have Auslander rigid dualizing complexes. A key step in the
proof of Theorem \ref{a0.5} is the following theorem
(see Theorem \ref{a7.1} for full details and proof).

\begin{thm} \label{a0.7}
Let $A$ be a noetherian connected graded $k$-algebra.
Suppose $t \in A$ is a nonzero homogeneous normal element of
positive degree. Then $A$ has a graded
Auslander balanced dualizing complex if and only if so does $A/(t)$.
\end{thm}

As noted on \cite{Zh} p.\ 399, the proof of \cite{SZ} Lemma 6.1(ii)
has a gap, and an alternative proof of the result (under extra
hypotheses) is given in \cite{Zh} Theorem 3.1. We now give
a complete proof of \cite{SZ} Lemma 6.1(ii) using Auslander dualizing
complexes (the proof is at end of Section 5).

\begin{prop} \label{a0.2}
Let $A$ be a noetherian locally finite $\mbb{N}$-graded $k$-algebra.
If $A$ is graded FBN, or $A$ has enough normal elements, then
$\opn{GKdim} M = \opn{Kdim} M \in {\mbb N}$
for every finitely generated left or right graded $A$-module $M$.
\end{prop}

Here are some other results we prove:
\begin{enumerate}
\item Gabber's Maximality Principle (Theorem \ref{a3.1}).
\item Existence of double-Ext spectral sequence (Proposition
\ref{a1.4}).
\item The existence of an Auslander rigid dualizing complex is
transferred to related algebras (Propositions \ref{a6.16} and
\ref{a6.18}, Corollaries \ref{a6.15}, \ref{a7.8} and \ref{a7.9}).
\end{enumerate}

The canonical dimension, denoted by $\opn{Cdim}$, is defined when $A$
has an Auslander dualizing complex (Definition \ref{a2.5}). It is an
exact finitely partitive dimension function (Theorem \ref{a2.19}).
Local duality implies that the canonical dimension is symmetric in the
graded case (Proposition \ref{a6.9}). Therefore if $A$ is a connected
graded algebra with an Auslander balanced dualizing complex, the
canonical dimension $\opn{Cdim}$ is exact, finitely partitive, and
symmetric on graded modules. Note that the Krull dimension, denoted by
$\opn{Kdim}$, is exact and finitely partitive, but it is unknown
whether it is symmetric. On the other hand Gelfand-Kirillov dimension,
denoted by $\opn{GKdim}$, is symmetric, but neither exact nor finitely
partitive in general. Hence the canonical dimension is the better
dimension function -- at least in the graded or filtered case.

The study of dualizing complexes over noncommutative rings presents
many interesting and subtle questions. We conclude the
introduction by mentioning two of them:

\begin{que}
Which (noetherian, affine) $k$-algebras have (rigid)
dualizing complexes?
\end{que}

\begin{que}
Is a rigid dualizing complex is always Auslander?
\end{que}

\bigskip \noindent {\em Acknowledgments.}\
We thank M. Artin and M. Van den Bergh for many
useful conversations and email correspondences.

\section{Dualizing complexes}

Let $k$ be a field and let $A$ be an associative $k$-algebra with $1$.
All $A$-modules will be by default left modules, and we
denote by $\cat{Mod} A$ the category of left $A$-modules.
Let $A^{\circ}$ be the opposite algebra, and let
$A^{\mrm{e}} := A \otimes A^{\circ}$ where $\otimes = \otimes_{k}$.
Thus an $A^{\mrm{e}}$-module $M$ is, in the conventional notation,
an $A$-$A$-bimodule ${}_{A}M_{A}$ central over $k$.
Most of our definitions and results have a left-right symmetry,
expressible by the exchange $A \leftrightarrow A^{\circ}$.
Since these symmetries are evident we shall usually not mention them.

Let $\msf{D}(\cat{Mod} A)$ be the derived category of $A$-modules,
and let $\msf{D}^{\star}(\cat{Mod} A)$, for
$\star = \mrm{b}, +, -$ or blank, be the full subcategories
of bounded, bounded below, bounded above, or unbounded complexes
respectively \cite{RD}.

Given another $k$-algebra $B$, the forgetful functor
$\cat{Mod}(A \otimes B^{\circ}) \ar \cat{Mod} A$
is exact, and so induces a functor
$\cat{D}^{\star}(\cat{Mod}(A \otimes B^{\circ})) \ar
\cat{D}^{\star}(\cat{Mod} A)$.
Now $A \otimes B^{\circ}$ is a projective $A$-module, so
any projective (resp.\ flat, injective) $(A \otimes B^{\circ})$-module
is projective (resp.\ flat, injective) over $A$.

Consider $k$-algebras $A, B, C$. For complexes
$M \in \cat{D}(\cat{Mod}(A \otimes B^{\circ}))$ and
$N \in \cat{D}(\cat{Mod}(A \otimes C^{\circ}))$,
with either
$M \in \msf{D}^{-}$ or $N \in \msf{D}^{+}$,
there is a derived functor
\[ \mrm{R} \opn{Hom}_{A}(M, N) \in
\cat{D}(\cat{Mod}(B \otimes C^{\circ})) . \]
It is calculated by replacing $M$ with an isomorphic complex
in $\cat{D}^{-}(\cat{Mod}(A \otimes B^{\circ}))$
which consists of projective modules over $A$; or
by replacing $N$ with an isomorphic complex
in $\cat{D}^{+}(\cat{Mod}(A \otimes C^{\circ}))$
which consists of injective modules over $A$.
For full details see \cite{RD} and \cite{Ye1}.
Note that for modules $M$ and $N$, viewed as complexes concentrated
in degree $0$, one has
\[ \mrm{H}^{q} \mrm{R} \opn{Hom}_{A}(M, N) =
\opn{Ext}^{q}_{A}(M, N) , \]
the latter being the usual Ext.

Because the forgetful functors
$\cat{Mod}(A \otimes B^{\circ}) \ar \cat{Mod} A$
etc.\ commute with \linebreak
$\mrm{R} \opn{Hom}_{A}(-,-)$
there is no need to mention them explicitly.

A complex $N \in \msf{D}^{+}(\cat{Mod} A)$ is said to have finite
injective dimension if there is an integer $q_{0}$ with
$\opn{Ext}^{q}_{A}(M, N) = 0$
for all $q > q_{0}$ and $M \in \cat{Mod} A$.

For the rest of this section $A$ denotes a left noetherian
$k$-algebra and $B$ denotes a right noetherian $k$-algebra.
(For instance we could take $A = B$ a two-sided noetherian algebra.)
Observe that the algebra $A \otimes B^{\circ}$ need not be left
noetherian.

The subcategory $\cat{Mod}_{\mrm{f}} A$
of finitely generated $A$-modules is abelian and closed
under extensions. Hence there is a full triangulated subcategory
$\msf{D}_{\mrm{f}}(\cat{Mod} A) \subset \msf{D}(\cat{Mod} A)$
consisting of all complexes with finitely generated cohomologies.

Dualizing complexes over commutative rings were introduced in
\cite{RD}. The noncommutative graded  version
first appeared in \cite{Ye1}, and we now give a slightly more general
version.

\begin{dfn}\label{a1.1}
Assume $A$ and $B$ are $k$-algebras, with $A$ left noetherian and $B$
right noetherian. A complex
$R \in \msf{D}^{\mrm{b}}(\cat{Mod}(A \otimes B^{\circ}))$
is called a {\em dualizing complex} if it satisfies the three
conditions below:
\begin{enumerate}
\rmitem{i} $R$ has finite injective dimension over $A$ and $B^{\circ}$.
\rmitem{ii} $R$ has finitely generated cohomology modules over
$A$ and $B^{\circ}$.
\rmitem{iii} The canonical morphisms
$B \ar \mrm{R} \opn{Hom}_{A}(R, R)$
in $\msf{D}(\cat{Mod} B^{\mrm{e}})$,
and
$A \ar \mrm{R} \opn{Hom}_{B^{\circ}}(R, R)$
in $\msf{D}(\cat{Mod} A^{\mrm{e}})$, are both isomorphisms.
\end{enumerate}
In case $A = B$, we shall say that $R$ is a dualizing complex
over $A$.
\end{dfn}

Condition (i) is equivalent to having an isomorphism
$R \cong I \in \cat{D}^{\mrm{b}}(\cat{Mod} A\otimes B^\circ)$,
where each $I^{q}$ is injective over $A$ and over $B^{\circ}$.

\begin{exa}
Suppose $A$ is commutative and $R$ is a dualizing complex in the
sense of \cite{RD}. If we consider $R$ as a complex of bimodules,
by identifying $A = A^{\circ}$, then $R$ is a dualizing complex
in the sense of the definition above. According to \cite{Ye3}, if
$\opn{Spec} A$ is connected, then any
dualizing complex $R'$ over $A$ is isomorphic to $R \otimes_{A} P[n]$,
where $P$ is an invertible bimodule (not necessarily central!) and
$n \in \mbb{Z}$.
\end{exa}

Some easy examples of dualizing complexes over noncommutative rings
are given in Example \ref{a2.2}.

The next proposition offers an explanation of the name ``dualizing
complex.'' The duality functors associated to $R$ are the
contravariant functors
\begin{alignat*}{3}
D  & := \mrm{R} \opn{Hom}_{A}(-, R) & & :
\cat{D}(\cat{Mod} A) & & \ar \cat{D}(\cat{Mod} B^{\circ}) \\
D^{\circ} & := \mrm{R} \opn{Hom}_{B^{\circ}}(-, R) & & :
\cat{D}(\cat{Mod} B^{\circ}) & & \ar \cat{D}(\cat{Mod} A) .
\end{alignat*}

\begin{prop} \label{a1.2}
Let $R \in \msf{D}(\cat{Mod}(A \otimes B^{\circ}))$
be a dualizing complex.
\begin{enumerate}
\item For any
$M \in \cat{D}_{\mrm{f}}(\cat{Mod} A)$
one has $D M \in \cat{D}_{\mrm{f}}(\cat{Mod} B^{\circ})$
and
$M \cong D^{\circ} D M$.
\item The functors $D$ and $D^{\circ}$ determine a duality, i.e.,\
an anti-equivalence, of triangulated categories between
$\msf{D}_{\mrm{f}}(\cat{Mod} A)$
and
$\msf{D}_{\mrm{f}}(\cat{Mod} B^{\circ})$,
restricting to a duality between
$\msf{D}^{\mrm{b}}_{\mrm{f}}(\cat{Mod} A)$
and
$\msf{D}^{\mrm{b}}_{\mrm{f}}(\cat{Mod} B^{\circ})$.
\end{enumerate}
\end{prop}

\begin{proof}
1.  This is slightly stronger than \cite{Ye1} Lemma 3.5.
By adjunction we get a functorial morphism
$M \ar D^{\circ} D M$.
Since the functor $D^{\circ} D$ is way
out in both directions and $D^{\circ} D A \cong A$ by assumption,
the claim follows from the reversed forms of
\cite{RD} Propositions I.7.1 and I.7.3.\\
2. Immediate from part 1, together with the fact that
$M \in \msf{D}^{\mrm{b}}(\cat{Mod} A)$ implies
$D M \in \msf{D}^{\mrm{b}}(\cat{Mod} B^{\circ})$.
\end{proof}

\begin{rem}\label{a1.3}
The noetherian hypothesis can be relaxed -- dualizing complexes can be
defined over any coherent algebra $A$ (see \cite{Ye1} or \cite{Mi1}).
The category of finitely generated $A$-modules is then replaced by the
category of coherent modules. Many definitions and results in our
paper hold for coherent algebras, as can be easily checked.

Perhaps one can even work over an arbitrary algebra, using the
category of coherent complexes, as defined by Illusie (see \cite{SGA6}
Expos\'{e} I).

Another direction to extend the theory is to allow $k$ to be any
commutative ring. In this case the derived category of bimodules
should be
$\msf{D}(\cat{DGMod}(A \otimes_{k}^{\mrm{L}} B^{\circ}))$,
where $A \otimes_{k}^{\mrm{L}} B^{\circ}$ is a differential graded
algebra. See \cite{Ye3} Remark 1.12.
\end{rem}

\begin{rem}
Miyachi proved a converse to Proposition \ref{a1.2}(2): if there are
contravariant triangle functors
$\cat{D}(\cat{Mod} A) \ar \cat{D}(\cat{Mod} B^{\circ})$
and
$\cat{D}(\cat{Mod} B^{\circ}) \ar \cat{D}(\cat{Mod} A)$
which send $\bigoplus$ to $\prod$, preserve $\cat{D}^{\mrm{b}}$
and induce duality on $\cat{D}^{\mrm{b}}_{\mrm{f}}$, then there is a
dualizing complex in
$\msf{D}(\cat{Mod}(A \otimes B^{\circ}))$
(see \cite{Mi2} Theorem 3.3).
\end{rem}

\begin{rem}
There are examples of algebras $A$ and $B$ where there is a dualizing
complex
$R \in \msf{D}(\cat{Mod}(A \otimes B^{\circ}))$,
but there is no dualizing complex in
$\msf{D}(\cat{Mod} A^{\mrm{e}})$;
cf.\ \cite{WZ}. The algebras $A$ and $B$ are necessarily not
derived Morita equivalent, since given a tilting complex
$T \in \msf{D}(\cat{Mod}(B \otimes A^{\circ}))$,
the complex
$R \otimes^{\mrm{L}}_{B} T \in \msf{D}(\cat{Mod} A^{\mrm{e}})$
would be dualizing (cf.\ \cite{Ye3}).
\end{rem}

There are Grothendieck spectral sequences for the isomorphism of
functors \linebreak
$\bsym{1}_{\msf{D}^{\mrm{b}}_{\mrm{f}}(\cat{Mod} A)}
\cong D^{\circ} D$
and
$\bsym{1}_{\msf{D}^{\mrm{b}}_{\mrm{f}}(\cat{Mod} B^{\circ})} \cong
D D^{\circ}$.
For modules they take this form:

\begin{prop} \label{a1.4}
Let
$R \in \msf{D}(\cat{Mod}(A \otimes B^{\circ}))$
be a dualizing complex. Then there are convergent
{\em double-Ext} spectral sequences
\begin{equation} \label{a1.5}
E_{2}^{p, q}: = \opn{Ext}^{p}_{B^{\circ}}(\opn{Ext}^{-q}_{A}
(M, R), R) \Rightarrow M
\end{equation}
for all $M \in \cat{Mod}_{\mrm{f}} A$, and
\begin{equation} \label{a1.6}
E_{2}^{p, q}: = \opn{Ext}^{p}_{A}(\opn{Ext}^{-q}_{B^\circ}
(N, R), R) \Rightarrow N
\end{equation}
for all $N \in \cat{Mod}_{\mrm{f}} B^\circ$.
\end{prop}

\begin{proof}
By symmetry it suffices to consider (\ref{a1.5}) only.
We can assume $R$ is a bounded complex of bimodules with each
$R^{q}$ an injective module over $A$ and $B^{\circ}$. Given a nonzero
finitely generated $A$-module $M$,
define the complex
\[ H := \opn{Hom}_{B^{\circ}}(\opn{Hom}_{A}(M, R), R) . \]
Then the adjunction homomorphism
$M \ar H$ is a quasi-isomorphism.
Pick a positive integer $d$ large enough so that
$R^{q} = 0$ if $\lvert q \rvert > d$. Consider the  decreasing
filtration on $H$ given by the subcomplexes
\[ F^{p} H := \opn{Hom}_{B^{\circ}}(\opn{Hom}_{A}(M, R), R^{\geq p})
. \]
Then $F$ is an exhaustive filtration, and it determines the
convergent spectral sequence (\ref{a1.5}).
\end{proof}

Given $M \in \msf{D}^{+}(\cat{Mod} A)$, there is a quasi-isomorphism
$M \ar I$ in  $\msf{D}^{+}(\cat{Mod} A)$, where each $I^{q}$
is injective and
$\opn{Ker}(I^{q} \ar I^{q + 1}) \subset I^{q}$
is essential.
Such $I$ is unique (up to a non-unique isomorphism), and it is
called the {\em minimal injective resolution} of $M$ (cf.\ \cite{Ye1}
Lemma 4.2). If $M$ has finite injective dimension then $I$ is
bounded.

The next two results are straightforward generalizations of
\cite{ASZ} Lemma 2.2 and Theorem 2.3, so the proofs are omitted.

\begin{lem} \label{a1.8}
Let
$R \in \msf{D}(\cat{Mod}(A \otimes B^{\circ}))$
be a dualizing complex, and let $I$ be the
minimal injective resolution of $R$ in $\msf{D}^{\mrm{b}}(\cat{Mod}
A)$. Let
$Z_{i} := \opn{Ker}(I^{i} \ar I^{i+1})$
and let $M$ be a finitely generated left $A$-module. Then there exist
$f_{1}, \ldots, f_{n} \in \opn{Hom}(M, Z_{i})$
such that for every $N \subset \bigcap_{j} \opn{Ker}(f_{j})$
the natural map
$\opn{Ext}^{i}_{A}(M, R) \ar \opn{Ext}^{i}_{A}(N, R)$
is zero; or equivalently, the natural map
$\opn{Ext}^{i}_{A}(M / N, R) \ar \opn{Ext}^{i}_{A}(M, R)$
is surjective.
\end{lem}

\begin{thm} \label{a1.9}
Let
$R \in \msf{D}(\cat{Mod}(A \otimes B^{\circ}))$
be a dualizing complex, let $I$ be the
minimal injective resolution of $R$ in
$\msf{D}^{\mrm{b}}(\cat{Mod} A)$, and let
$Z_{i} := \opn{Ker}(I^{i} \ar I^{i+1})$.
Then:
\begin{enumerate}
\item For every nonzero $A$-module $M$ there is a nonzero
submodule $N \subset M$ which embeds in some $Z_{i}$.
\item Every indecomposable  injective $A$-module appears in $I$.
\end{enumerate}
\end{thm}

We conclude this section with a discussion of dualizing complexes
in $\msf{D}(\cat{Mod}(A \otimes B^{\circ}))$
when $A$ is commutative. For a prime ideal
$\mfrak{p} \subset A$ let $J_{A}(\mfrak{p})$ be an injective hull of
$A / \mfrak{p}$. Let us recall a result of Grothendieck.

\begin{prop}[\cite{RD} Proposition V.7.3] \label{b1.10}
Suppose $A$ is a commutative noetherian ring and
$R \in \msf{D}^{\mrm{b}}_{\mrm{f}}(\cat{Mod} A)$
is a \tup{(}central\tup{)} dualizing complex.
Let $I$ be the minimal injective resolution of $R$ in
$\cat{Mod} A$.
\begin{enumerate}
\item There is a function $d : \opn{Spec} A \ar \mbb{Z}$ such that
\[ I^{q} \cong \bigoplus_{d(\mfrak{p}) = q} J_{A}(\mfrak{p}) . \]
\item If $\mfrak{p} \subset \mfrak{q}$ are primes and
$\mfrak{q} / \mfrak{p} \subset A / \mfrak{p}$ has height $1$,
then $d(\mfrak{p}) = d(\mfrak{q}) - 1$.
\end{enumerate}
\end{prop}

The function $d$ is called a codimension function in \cite{RD}.
In our case we get:

\begin{thm} \label{a1.10}
Suppose $A$ is a commutative noetherian $k$-algebra,
$B$ is a right noetherian $k$-algebra,
and $R \in \msf{D}(\cat{Mod}(A \otimes B^{\circ}))$ is a dualizing
complex.
Let $I$ be the minimal injective resolution of $R$ in
$\cat{Mod} A$. Then:
\begin{enumerate}
\item There are functions
$d, r : \opn{Spec} A \ar \mbb{Z}$,
with $r \geq 1$ and constant on connected components of
$\opn{Spec} A$, s.t.\
\[ I^{q} \cong \bigoplus_{d(\mfrak{p}) = q}
J_{A}(\mfrak{p})^{r(\mfrak{p})} . \]
\item If $\mfrak{p} \subset \mfrak{q}$ are primes of $A$ and
$\mfrak{q} / \mfrak{p} \subset A / \mfrak{p}$ has height $1$,
then $d(\mfrak{p}) = d(\mfrak{q}) - 1$.
\item $A$ is catenary, and \tup{(}if $A \neq 0$\tup{)} its Krull
dimension is
\[ \opn{Kdim} A \leq \max \{ d(\mfrak{p}) \} - \min \{ d(\mfrak{p}) \}
< \infty . \]
\item $B$ is an Azumaya $A$-algebra.
\end{enumerate}
\end{thm}

The proof of the theorem is after the next lemma.

\begin{lem} \label{c2.1}
Assume in addition that $A$ is local. Then there is an integer $d$
such that
$\opn{Ext}^{q}_{A}(M, R) = 0$
and
$\opn{Ext}^{q}_{B^{\circ}}(N, R) = 0$
for all $q \neq d$, all finite length $A$-modules $M$ and all
finite length $B^{\circ}$-modules $N$.
\end{lem}

\begin{proof}
According to \cite{Ye3} Proposition 5.4 (which works even when
$A \neq B$; cf.\ ibid.\ Proposition 2.5), left and right
multiplications on $R$ induce ring isomorphisms
$A \cong
\opn{End}_{\cat{D}(\cat{Mod} (A \otimes B^{\circ}))} (R) \cong
\opn{Z}(B)$.
Moreover since $A$ is noetherian and
$B^{\circ} \cong \opn{Ext}^{0}_{A}(R, R)$
we see that $B$ is a finite $A$-algebra.
If $N$ is an $A$-central $(A \otimes B^{\circ})$-module,
then $\opn{Ext}^{q}_{B^{\circ}}(N, R)$
is a central $A$-bimodule.

Denote by $K$ the residue field of $A$.
Let
$p_{0} := \opn{min} \{ p \mid \opn{Ext}^{p}_{A}(K, R) \neq 0 \}$
and
$p_{1} := \opn{max} \{ p \mid \opn{Ext}^{p}_{A}(K, R) \neq 0 \}$.
By induction on length we see that for every finite length
$A$-module $M$,
$p_{0} = \opn{min} \{ p \mid \opn{Ext}^{p}_{A}(M, R) \neq 0 \}$
and
$p_{1} = \opn{max} \{ p \mid \opn{Ext}^{p}_{A}(M, R) \neq 0 \}$.

Now take a nonzero finite length $B^{\circ}$-module $N$. Since
we can view $N$ as a central $A$-bimodule, it follows that
$\opn{Ext}^{q}_{B^{\circ}}(N, R)$
is also a central $A$-bimodule for every $q$; and so it
has finite length. Define
$q_{0} := \opn{min} \{ q \mid \opn{Ext}^{q}_{B^{\circ}}
(N, R) \neq 0 \}$
and
$q_{1} := \opn{max} \{ q \mid \opn{Ext}^{q}_{B^{\circ}}
(N, R) \neq 0 \}$.
In the $E_{2}$-page of the spectral sequence of Proposition
\ref{a1.4} we have nonzero terms
$E_{2}^{p_{1}, -q_{0}} = \opn{Ext}^{p_{1}}_{A}
(\opn{Ext}^{q_{0}}_{B^{\circ}}(N, R))$
and
$E_{2}^{p_{0}, -q_{1}} = \opn{Ext}^{p_{0}}_{A}
(\opn{Ext}^{q_{1}}_{B^{\circ}}(N, R))$,
that appear in the right-top and left-bottom corners
respectively. The convergence of the spectral sequence forces
$p_{1} = q_{0}$ and $p_{0} = q_{1}$.
Therefore we get
$d := p_{1} = q_{0} = q_{1} = p_{0}$.
\end{proof}

\begin{proof}[Proof of Theorem \tup{\ref{a1.10}}]
1. As observed in the proof of the lemma, $B$ is a finite
$A$-algebra. Take a prime $\mfrak{p} \subset A$, and define
$B_{\mfrak{p}} := B \otimes_{A} A_{\mfrak{p}}$
and
$R_{\mfrak{p}} := A_{\mfrak{p}} \otimes_{A} R \otimes_{B}
B_{\mfrak{p}} \in
\cat{D}^{\mrm{b}}(\cat{Mod}(A_{\mfrak{p}} \otimes
B^{\circ}_{\mfrak{p}}))$.
We claim that $R_{\mfrak{p}}$ is a dualizing complex.

First note that the cohomologies $\opn{H}^{q} R$ are central
$A$-bimodules. Since $A \ar A_{\mfrak{p}}$ is flat and
$A_{\mfrak{p}} \otimes_{A} A_{\mfrak{p}} = A_{\mfrak{p}}$,
there are isomorphisms
$\opn{H}^{q} R_{\mfrak{p}} \cong
A_{\mfrak{p}} \otimes_{A} \opn{H}^{q} R \cong
\opn{H}^{q} R \otimes_{B} B_{\mfrak{p}}$.
Therefore
$R_{\mfrak{p}} \cong A_{\mfrak{p}} \otimes_{A} R
\cong R \otimes_{B} B_{\mfrak{p}}$
in $\cat{D}^{\mrm{b}}(\cat{Mod}(A \otimes B^{\circ}))$.
It is easy to see that the cohomology bimodules $\opn{H}^{q} R$
are finitely generated on both sides,
$\opn{R} \opn{Hom}_{A_{\mfrak{p}}}(R_{\mfrak{p}}, R_{\mfrak{p}})
\cong B^{\circ}_{\mfrak{p}}$
and
$\opn{R} \opn{Hom}_{B^{\circ}_{\mfrak{p}}}
(R_{\mfrak{p}}, R_{\mfrak{p}}) \cong A_{\mfrak{p}}$.

In order to verify that $R_{\mfrak{p}}$ has finite injective 
dimension over $B^{\circ}_{\mfrak{p}}$
it suffices to show that
$\opn{Ext}^{q}_{B^{\circ}_{\mfrak{p}}}(N, R_{\mfrak{p}})$
vanishes for all finitely generated 
$B^{\circ}_{\mfrak{p}}$-modules $N$, for large $q$.
Now we can write
$N \cong N' \otimes_{B} B_{\mfrak{p}}$
for some finitely generated $B^{\circ}$-module $N'$,
and then
$\opn{Ext}^{q}_{B^{\circ}_{\mfrak{p}}}(N, R_{\mfrak{p}}) \cong
\opn{Ext}^{q}_{B^{\circ}}(N', R) \otimes_{A} A_{\mfrak{p}}$.
Likewise for the injective dimension over$A_{\mfrak{p}}$. 
So indeed $R_{\mfrak{p}}$ is dualizing.

Since the multiplicity of the indecomposable injective
$J_{A}(\mfrak{p})$ in $I$ is measured by
$\opn{Ext}^{q}_{A_{\mfrak{p}}}(k(\mfrak{p}), R_{\mfrak{p}})$,
where $k(\mfrak{p})$ is the residue field, the lemma says that
$J_{A}(\mfrak{p})$  occurs in the complex $I$ in only one degree,
say $d(\mfrak{p})$. The fact that the multiplicity $r(\mfrak{p})$
is locally constant will be proved in part 4 below.

\medskip \noindent
2. Choose $a \in \mfrak{q} - \mfrak{p}$, so in the exact sequence
\[ 0 \ar (A / \mfrak{p})_{\mfrak{q}} \xrightarrow{a}
(A / \mfrak{p})_{\mfrak{q}} \ar M \ar 0 \]
the $A_{\mfrak{q}}$-module $M$ has finite length. Applying
$\opn{Ext}^{q}_{A}(-, R)$ to this sequence we obtain an exact
sequence of finitely generated $A_{\mfrak{q}}$-modules
\[ \opn{Ext}^{q}_{A}((A / \mfrak{p})_{\mfrak{q}}, R)
\xrightarrow{a}
\opn{Ext}^{q}_{A}((A / \mfrak{p})_{\mfrak{q}}, R) \ar
\opn{Ext}^{q + 1}_{A}(M, R) \]
for each $q$. By Nakayama's Lemma and part 1 we find that
$\opn{Ext}^{q}_{A}((A / \mfrak{p})_{\mfrak{q}}, R) = 0$
unless $q + 1 = d(\mfrak{q})$. Hence
$d(\mfrak{p}) = d(\mfrak{q}) - 1$
as claimed.

\medskip \noindent
3. This follows trivially from 2.

\medskip \noindent
4. Pick a prime ideal $\mfrak{p} \subset A$. Since $R_{\mfrak{p}}$
is a dualizing complex in
$\cat{D}(\cat{Mod}(A_{\mfrak{p}} \otimes B_{\mfrak{p}}^{\circ}))$,
the lemma tells us that that the functors
$\opn{Ext}^{d(\mfrak{p})}_{A_{\mfrak{p}}}(-, R_{\mfrak{p}})$
and
$\opn{Ext}^{d(\mfrak{p})}_{B_{\mfrak{p}}^{\circ}}(-, R_{\mfrak{p}})$
are a duality between the categories of finite length modules over
$A_{\mfrak{p}}$ and $B_{\mfrak{p}}^{\circ}$. Furthermore since
this duality is $A_{\mfrak{p}}$-linear, it restricts to a duality
between
$\cat{Mod}_{\mrm{f}} A_{\mfrak{p}} / \mfrak{p}_{\mfrak{p}}^{n}$
and
$\cat{Mod}_{\mrm{f}}(B_{\mfrak{p}} / \mfrak{p}_{\mfrak{p}}^{n}
B_{\mfrak{p}})^{\circ}$
for every $n \geq 1$.
Since
$\cat{Mod}_{\mrm{f}} A_{\mfrak{p}} / \mfrak{p}_{\mfrak{p}}^{n}$
has an auto-duality, namely
$\opn{Hom}_{A_{\mfrak{p}}}(-, J_{A}(\mfrak{p}))$,
it follows that there is an equivalence between
$\cat{Mod}_{\mrm{f}} A_{\mfrak{p}} / \mfrak{p}_{\mfrak{p}}^{n}$
and
$\cat{Mod}_{\mrm{f}} B_{\mfrak{p}} / \mfrak{p}_{\mfrak{p}}^{n}
B_{\mfrak{p}}$.
Morita Theory says there is an isomorphism
$\phi_{n} : B_{\mfrak{p}} / \mfrak{p}_{\mfrak{p}}^{n}
B_{\mfrak{p}} \iso
\opn{M}_{r(n)}(A_{\mfrak{p}} / \mfrak{p}_{\mfrak{p}}^{n})$
for some number $r(n)$.

Since the isomorphisms $\phi_{n}$ all arise from the same
equivalence they can be made compatible, and in particular we have
$r(n) = r(\mfrak{p})$ for all $n$. In the inverse limit we get
$\hat{A}_{\mfrak{p}} \otimes_{A} B \cong
\opn{M}_{r(\mfrak{p})}(\hat{A}_{\mfrak{p}})$
as $\hat{A}_{\mfrak{p}}$-algebras.

Since $\mfrak{p}$ was an arbitrary prime ideal of $A$, we conclude
that the multiplication map
$B \otimes_{A} B^{\circ} \ar \opn{End}_{A}(B)$
is bijective. Therefore $B$ is Azumaya over $A$. As a projective
$A$-module, the rank of $B$ at a prime $\mfrak{p}$ is precisely
$r(\mfrak{p})^{2}$, so the function $r$ is locally constant.
\end{proof}

\begin{rem}
A complex such as $I$ in Theorem \ref{b1.10} is called a residual
complex. It actually depends functorially on $R$: $I = \mrm{E} R$,
where $\mrm{E}$ is the Cousin functor. Noncommutative variants of the
Cousin functor are studied in \cite{Ye2} and \cite{YZ2}. In
particular one can show that the complex $I$ can
be made into a complex of bimodules, where on the right it is the
minimal injective resolution of $R$ in $\cat{Mod} B^{\circ}$.
It follows that $R$ is an Auslander dualizing complex, as defined in
Section 2. In \cite{YZ2} we show that a more complicated  version of
Theorem \ref{a1.10} holds when $A = B$ is a PI algebra.
\end{rem}

\begin{exa} \label{a1.12}
Let $A$ be a noetherian commutative regular ring of infinite Krull
dimension -- see Nagata's example \cite{Na} Appendix, p.\ 203
Example 1.
The complex $R = A$ is a pointwise dualizing complex (see \cite{RD}
Section V.8), and
$\mrm{R} \opn{Hom}_{A}(-, A): \msf{D}_{\mrm{f}}
(\cat{Mod} A) \ar \msf{D}_{\mrm{f}}(\cat{Mod} A)$
is a duality. However by Theorem \ref{a1.10}, there is no
dualizing complex in
$\msf{D}^{\mrm{b}}(\cat{Mod}(A \otimes B^{\circ}))$
for any right noetherian $k$-algebra $B$.
\end{exa}

\begin{exa}
Let $A$ be the non-catenary noetherian commutative local ring of
\cite{Na} Appendix, p.\ 203 Example 2. Then again there is no
dualizing complex in
$\msf{D}^{\mrm{b}}(\cat{Mod}(A \otimes B^{\circ}))$
for any right noetherian $k$-algebra $B$.
\end{exa}

\section{Auslander Dualizing Complexes}

The basic ideas in this section already appear in \cite{Ye2} Section
1 (which treats graded algebras).

Assume
$R \in \cat{D}^{\mrm{b}}(\cat{Mod} (A \otimes B^{\circ}))$
is a dualizing complex.
Let $M$ be a finitely generated $A$-module. The
{\em grade of $M$ with respect to $R$} is
\[ j_{R; A}(M) := \inf \{ q \mid \opn{Ext}^{q}_{A}(M, R) \neq 0 \}
\in \mbb{Z} \cup \{\infty\} . \]
Similarly define $j_{R; B^{\circ}}$ for a $B^{\circ}$-module.

We are ready to define the notion appearing in the title of the paper.

\begin{dfn} \label{a2.1}
Let $A$ and $B$ be $k$-algebras, with $A$ left noetherian and
$B$ right noetherian, and let
$R \in \msf{D}(\cat{Mod}(A \otimes B^{\circ}))$
be a dualizing complex.
We say that $R$ has the {\em Auslander property}, or that $R$ is
an {\em Auslander dualizing complex}, if
\begin{enumerate}
\rmitem{i} for every finitely generated $A$-module $M$,
integer $q$ and $B^{\circ}$-submodule
$N \subset \opn{Ext}^{q}_{A}(M, R)$,
one has $j_{R; B^{\circ}}(N) \geq q$.
\rmitem{ii} the same holds after exchanging  $A$ and $B^{\circ}$.
\end{enumerate}
\end{dfn}

Note that the role of the algebras $A$ and $B^{\circ}$ is symmetric.
Also note that if $R$ is an Auslander dualizing complex, then any
shift $R[n]$ is also an Auslander dualizing complex (the shift cancels
out in the double dual).

\begin{exa}
Let $A$ be a commutative $k$-algebra and $R$ a central dualizing
complex over it. From Proposition \ref{b1.10} it is clear that $R$ is
Auslander. For a prime ideal $\mfrak{p}$ one has
$j_{R; A}(A / \mfrak{p}) = d(\mfrak{p})$.
\end{exa}

\begin{exa} \label{a2.2}
(a) If $A$ is Gorenstein, i.e.\ the bimodule $A$ has
finite injective dimension on both sides, then $R = A$ is a dualizing
complex over $A$. If $A$ is an Auslander-Gorenstein ring in
the sense of \cite{Bj}, then $R=A$ is an Auslander dualizing
complex.

\noindent
(b) If $A$ is a finite $k$-algebra, i.e.\
$\opn{rank}_{k} A < \infty$, then the bimodule
$A^{*} := \opn{Hom}_{k}(A, k)$ is injective on both sides,
and $R = A^{*}$ is an Auslander dualizing complex over $A$.
Clearly $j_{R; A}(M) = 0$ for all $M$.

\noindent
(c) Let $A$ be the matrix algebra
$\left( \begin{smallmatrix} k & V \\ 0 & k \end{smallmatrix}
\right)$
where $V$ is a finite rank $k$-module. Since $A$ is hereditary it
is also Gorenstein, so both $A$ and $A^{*}$ are dualizing complexes,
and if $V \neq 0$ they are non-isomorphic. According to \cite{ASZ}
Example 5.4, the dualizing complex $A$ is Auslander iff
$\opn{rank}_{k} V \leq 1$. See also \cite{Ye3} Section 3.
\end{exa}

More examples of algebras with Auslander dualizing complexes are in
Examples \ref{a8.26}-\ref{exa8.1}.

The next definition is taken from \cite{MR} Section 6.8.4 (with a
slight modification -- we allow negative dimensions).

\begin{dfn} \label{a2.3}
An {\em exact dimension function} is a function $\opn{dim}$ which
assigns to each module $M \in \cat{Mod}_{\mrm{f}} A$
a value $\opn{dim} M$ in an ordered set containing $-\infty$,
$\mbb{R}$ and some infinite ordinals, and satisfies the following
axioms:
\begin{enumerate}
\rmitem{i} $\opn{dim} 0 = -\infty$.
\rmitem{ii} For every short exact sequence
$0 \ar M' \ar M \ar M'' \ar 0$ one has
$\opn{dim} M = \max \{ \opn{dim} M', \opn{dim} M'' \}$.
\rmitem{iii} If $\mfrak{p} M = 0$ for some prime ideal $\mfrak{p}$,
and $M$ is a torsion $A / \mfrak{p}$-module, then
$\opn{dim}_{R} M \leq \opn{dim}_{R} A / \mfrak{p} - 1$.
\end{enumerate}
\end{dfn}

Familiar examples of exact dimension functions are the Krull dimension
$\opn{Kdim}$ and (sometimes) the Gelfand-Kirillov dimension
$\opn{GKdim}$.

\begin{lem} \label{a2.17}
A function $\opn{dim}$ defined on $\cat{Mod}_{\mrm{f}} A$ and
satisfying axioms \tup{(i)-(ii)}
\tup{(}resp.\ axioms \tup{(i)-(iii))}
extends uniquely to a function on $\cat{Mod} A$,
satisfying axioms \tup{(i)-(ii)}
\tup{(}resp.\ axioms \tup{(i)-(iii))} and the axiom
\begin{enumerate}
\rmitem{iv} $\opn{dim} M = \sup \{ \opn{dim} M' \mid M' \subset M
\text{ finitely generated} \}$.
\end{enumerate}
\end{lem}

The proof of the lemma is standard.

Usually we will have a pair of dimension functions, one on
$\cat{Mod} A$ and the other on $\cat{Mod} B^{\circ}$;
when necessary we shall
distinguish between them by writing $\opn{dim}_{A}$ and
$\opn{dim}_{B^{\circ}}$ respectively.

\begin{dfn}
An exact dimension function $\opn{dim}$ is called {\em finitely
partitive} if given a finitely generated
module $M$ there is a number $l_{0}$, such that for every chain of
submodules
$M = M_{0} \supsetneqq M_{1} \cdots \supsetneqq M_{l}$
with
$\opn{dim} M_{i} / M_{i + 1} = \opn{dim} M$
one has $l \leq l_{0}$.
\end{dfn}

In the next two definitions (taken from \cite{ASZ} and \cite{Ye2})
$\opn{dim}$ denotes a function on $\cat{Mod} A$ satisfying axioms
(i), (ii) and (iv).

\begin{dfn}
\begin{enumerate}
\item A module $M$ is called {\em $\opn{dim}$-pure} if
$\opn{dim} M' = \opn{dim} M$ for every nonzero submodule
$M' \subset M$.
\item A module $M$ is called {\em $\opn{dim}$-essentially pure}
if $M$ contains an essential submodule which is pure.
\item A module $M$ is called {\em $\opn{dim}$-critical} if
$M \neq 0$, and $\opn{dim} M / M' < \opn{dim} M$ for every
$0 \neq M' \subsetneqq M$.
\end{enumerate}
\end{dfn}

\begin{dfn} \label{a2.4}
\begin{enumerate}
\item Let $\msf{M}_{q}(\opn{dim})$ be the full subcategory of
$\cat{Mod} A$ consisting of modules $M$ with $\opn{dim} M \leq q$,
and let
$\msf{M}_{q, \mrm{f}}(\opn{dim}) :=
\msf{M}_{q}(\opn{dim}) \cap \msf{Mod}_{\mrm{f}} A$.
\item Given a module $M$ let
$\Gamma_{\msf{M}_{q}(\opn{dim})} M \subset M$
be the largest submodule $M' \subset M$ such that
$\opn{dim} M' \leq q$.
\end{enumerate}
\end{dfn}

Since $\msf{M}_{q}(\opn{dim})$ is a localizing subcategory of
$\msf{Mod} A$ the submodule $\Gamma_{\msf{M}_{q}(\opn{dim})} M$
is well defined (and in fact $\Gamma_{\msf{M}_{q}(\opn{dim})}$ is a
left exact idempotent functor).
The corresponding subcategories of $\cat{Mod} B^{\circ}$ shall be
denoted by $\msf{M}^{\circ}_{q}(\opn{dim})$ and
$\msf{M}^{\circ}_{q, \mrm{f}}(\opn{dim})$.

\begin{dfn}\label{a2.5}
Let $M$ be a finitely generated $A$-module.
The {\em canonical dimension of $M$ with respect to $R$} is
\[ \opn{Cdim}_{R; A} M := - j_{R; A}(M) \in \mbb{Z} \cup
\{ -\infty \}  . \]
Likewise define $\opn{Cdim}_{R; B^{\circ}}$.
\end{dfn}

The canonical dimension will not be an exact dimension function in
general. However we have the following theorem, which generalizes
results of Bj\"ork and Levasseur (the graded case was proved in
\cite{Ye2}).

\begin{thm} \label{a2.19}
If $R$ is an Auslander dualizing complex then
$\opn{Cdim}_{R; A}$ is a finitely partitive exact dimension function.
\end{thm}

The proof of this theorem appears later on in the section. The key
step is:

\begin{lem} \label{a2.21}
Let $0 \ar M' \ar M' \ar M'' \ar 0$ be a short exact sequence
 of finitely generated $A$-modules. Then
 \[ j_{R; A}(M) = \inf \{ j_{R; A}(M'), j_{R; A}(M'') \} . \]
\end{lem}

\begin{proof}
The proof goes along the lines of the proofs in \cite{Bj} and
\cite{Lev}.
By Proposition \ref{a1.4} we have a convergent spectral sequence
\begin{equation} \label{a2.16}
E_{2}^{p, q}: = \opn{Ext}^{p}_{B^{\circ}}(\opn{Ext}^{-q}_{A}
(M, R), R) \Rightarrow M ,
\end{equation}
so there is a corresponding filtration (called the b-filtration in
\cite{Lev} Theorem 2.2)
\[ M = F^{-d} M \supset F^{-d+1} M \supset \cdots \supset
F^{d+1} M = 0. \]

The Auslander condition tells us that
$E_{2}^{p, q} = 0$ if $p < -q$. So the spectral sequence lives
in a bounded region of the $(p,q)$ plane:
$p \geq -q$, $q \leq -j_{A; R}(M)$ and $p \leq d$
(see Figure 1).
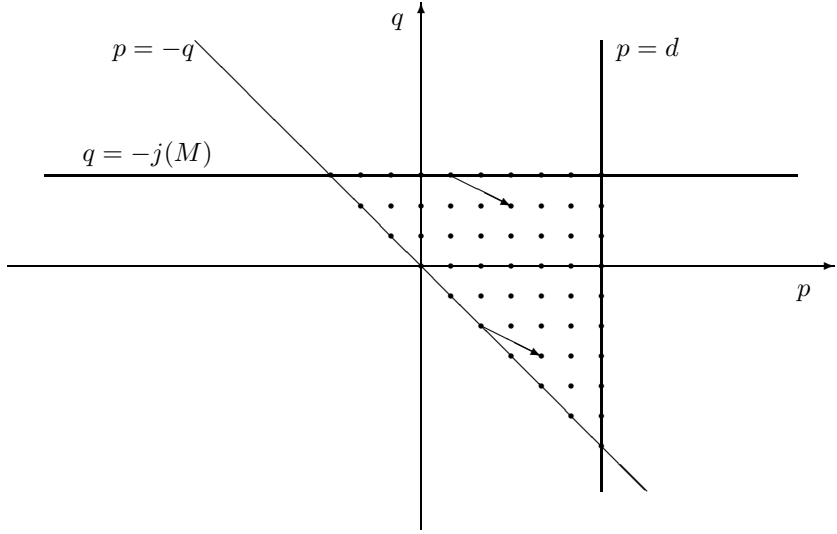
\begin{figure} \label{fig1}
\setlength{\unitlength}{1mm}
\begin{picture}(110,70)(-55,-35)

\put(-55,0){\vector(1,0){110}}
\put(50,-4){$p$}

\put(0,-35){\vector(0,1){70}}
\put(-4,32){$q$}

\put(-30,30){\line(1,-1){60}}
\put(-41,28){$p = -q$}
\put(-50,12){\line(1,0){100}}
\put(-45,14){$q = -j(M)$}
\put(24,-30){\line(0,1){60}}
\put(26,28){$p = d$}

\multiput(-12,12)(4,0){10}{\circle*{.8}}
\multiput(-8,8)(4,0){9}{\circle*{.8}}
\multiput(-4,4)(4,0){8}{\circle*{.8}}
\multiput(0,0)(4,0){7}{\circle*{.8}}
\multiput(4,-4)(4,0){6}{\circle*{.8}}
\multiput(8,-8)(4,0){5}{\circle*{.8}}
\multiput(12,-12)(4,0){4}{\circle*{.8}}
\multiput(16,-16)(4,0){3}{\circle*{.8}}
\multiput(20,-20)(4,0){2}{\circle*{.8}}
\multiput(24,-24)(4,0){1}{\circle*{.8}}

\put(4,12){\vector(2,-1){8}}
\put(8,-8){\vector(2,-1){8}}

\end{picture}
\caption{The $E_{2}$ term of the spectral sequence (\ref{a2.16})
in the $(p, q)$ plane}
\end{figure}
The coboundary operator of $E_{r}$ has bidegree $(r, 1 - r)$
and $r \geq 2$. We conclude
that for every $\lvert p \rvert \leq d$ there is an exact
sequence of $A$-modules
\begin{equation} \label{a2.7}
0 \ar \frac{F^{p} M}{F^{p+1} M} \ar E_{2}^{p, -p}
\ar Q^{p} \ar 0
\end{equation}
with $Q^{p}$ a subquotient of
$\bigoplus_{r \geq 2} E_{2}^{p + r, -p + (1 - r)}$
(cf.\ \cite{Bj} Theorem 1.3 and \cite{Lev} Theorem 2.2).

By the Auslander property it then follows that
$j_{R;A}(F^{p} M / F^{p+1} M) \geq p$ for all $M$ and $p$.
Just as in \cite{Bj} Proposition 1.6, one uses descending induction
on $p$, starting at $p = d + 1$, to prove that
$j_{R;A}(F^{p} M) \geq p$ for all $p$.
This implies that
\[ j_{R;B^\circ}(\opn{Ext}^{j_{R;A}(M)}(M, R)) = j_{R;A}(M) . \]
Now continue exactly like in \cite{Bj} Proposition 1.8.
\end{proof}

We conclude that $\opn{Cdim}_{R; A}$ verifies axiom (ii). Axiom (i)
holds trivially. By symmetry $\opn{Cdim}_{R; B^{\circ}}$ also
verifies axioms (i)-(ii).

\begin{thm} \label{a2.6}
Suppose $R \in \msf{D}(\cat{Mod}(A \otimes B^{\circ}))$
is an Auslander dualizing complex.
Let $M \in \cat{Mod}_{\mrm{f}} A$ be nonzero
and $\opn{Cdim}_{R} M = n$. Then:
\begin{enumerate}
\item $\opn{Ext}^{-n}_{A}(M, R)$ is $\opn{Cdim}_{R}$-pure
of dimension $n$.
\item For each $p$,
$\opn{Ext}^{-p}_{B^{\circ}}(\opn{Ext}^{-p}_{A}(M, R), R)$
is $\opn{Cdim}_{R}$-pure of dimension $p$, or is $0$.
\item For each $p$ there is an exact sequence
\[ 0 \ar \Gamma_{\msf{M}_{p - 1}} M \ar 
\Gamma_{\msf{M}_{p}} M \ar
\opn{Ext}^{-p}_{B^{\circ}}(\opn{Ext}^{-p}_{A}(M, R), R)
\ar Q^{p} \ar 0 , \]
functorial in $M$, where $\opn{Cdim}_{R} Q^{p} \leq p - 2$
and
$\msf{M}_{p} = \msf{M}_{p}(\opn{Cdim}_{R})$.
\end{enumerate}
\end{thm}

\begin{proof}
1. Because the line
$q = -j_{R;A}(M)$ is on the boundary of the region of
support of the spectral sequence (\ref{a2.16}),
and the coboundary operator of $E_{r}$ has bidegree $(r, 1 - r)$,
$r \geq 2$, we see that for this value of $q$ there is a
bounded filtration
$E_{2}^{p, q} \supset E_{3}^{p, q} \supset \cdots$,
with
$E_{r}^{p, q} / E_{r + 1}^{p, q}$ a subquotient of
$E_{r}^{p + r, q + (1 - r)}$.
Now the abutment of the spectral sequence
is concentrated on the line $p = -q$ of the $(p, q)$-plane,
so $E_{r}^{p, q} = 0$ for $p > -q$ and $r \gg 0$.
By Lemma \ref{a2.21} we conclude that for $q = -j_{R; A}(M)$ 
and $p > -q$,
\[ j_{R; A}(\opn{Ext}^{p}_{B^{\circ}}(\opn{Ext}^{-q}_{A}(M, R), R)) =
j_{R; A}(E_{2}^{p, q}) \geq p + 2 \]
(cf.\ \cite{Bj} formula (1.10)).
Just like in \cite{Bj} Proposition 1.11 it follows that \linebreak
$\opn{Ext}^{p}_{A}(\opn{Ext}^{p}_{B^{\circ}}(N, R), R) = 0$
for $p > j_{R; A}(M)$ and
$N = \opn{Ext}^{j_{R;A}(M)}_{A}(M, R)$.
So by \cite{Bj} Proposition 1.9 we conclude that $N$ is pure.

\medskip \noindent 2.
Take $N := \opn{Ext}^{-p}_{A}(M, R)$. Then
$\opn{Cdim}_{R;B^\circ} N \leq p$ and part 1 applies.

\medskip \noindent 3.
By part 2, the sequence (\ref{a2.7}) and induction on $p$ 
we see that
$\Gamma_{\msf{M}_{p}} M = F^{-p} M$.
\end{proof}

For an integer $q$ let
$\msf{M}_{q} := \msf{M}_{q}(\opn{Cdim}_{R}) \subset \cat{Mod} A$
be the localizing subcategory from Definition \ref{a2.4}.
The filtration by dimension of support $\{ \msf{M}_{q} \}$
of $\cat{Mod} A$ is called the {\em niveau filtration} in commutative
algebraic geometry. For each $q$ the quotient category
$\msf{M}_{q} / \msf{M}_{q - 1}$ is a locally noetherian abelian
category, and the full subcategory
$\msf{M}_{q, \mrm{f}} / \msf{M}_{q - 1, \mrm{f}}$ is noetherian
(see \cite{ASZ} Lemma 1.1). By symmetry we have
corresponding localizing subcategories
$\msf{M}^\circ_{q, \mrm{f}} \subset \msf{M}^\circ_{q}
\subset \cat{Mod} B^\circ$.

Recall from \cite{ASZ} Section 1 that two abelian categories $\cat{C}$
and $\cat{D}$ are said to be {\em dual} if they are
anti-equivalent, i.e.\ if $\cat{C}$ is equivalent to the opposite
category $\cat{D}^{\circ}$. Two categories $\cat{C}$ and
$\cat{D}$ are said to be in {\em step duality} if there are
filtrations by dense abelian subcategories
\[ 0 = \cat{C}_{n_{0} - 1} \subset \cat{C}_{n_{0}} \subset \cdots
\subset \cat{C}_{n_{1}} = \cat{C}
\quad {\text{and}} \quad
0 = \cat{D}_{n_{0} - 1} \subset \cat{D}_{n_{0}} \subset \cdots
\subset \cat{D}_{n_{1}} = \cat{D} \]
such that the quotient categories $\cat{C}_{i} / \cat{C}_{i - 1}$
and $\cat{D}_{i} / \cat{D}_{i - 1}$ are dual for all
$i = n_{0}, \ldots, n_{1}$.
Now Theorem \ref{a0.3} is a special case of:

\begin{thm} \label{a2.9}
Suppose $R \in \msf{D}(\cat{Mod}(A \otimes B^{\circ}))$
is an Auslander dualizing complex.
Then $\cat{Mod}_{\mrm{f}} A$ and $\cat{Mod}_{\mrm{f}} B^{\circ}$
are in step duality. More precisely, for every $q$ the functors
$\opn{Ext}^{q}_{A}(-, R)$ and
$\opn{Ext}^{q}_{B^{\circ}}(-, R)$
induce a duality between
the quotient categories
$\msf{M}_{q, \mrm{f}} / \msf{M}_{q - 1, \mrm{f}}$
and
$\msf{M}^{\circ}_{q, \mrm{f}} / \msf{M}^{\circ}_{q - 1, \mrm{f}}$.
\end{thm}

\begin{proof} Use Theorem \ref{a2.6} and the proof of \cite{ASZ}
Theorem 1.2.
\end{proof}

\begin{cor} \label{a2.10}
For each $q$ the category
$\msf{M}_{q, \mrm{f}} / \msf{M}_{q - 1, \mrm{f}}$ is artinian,
i.e.\ every object has finite length.
\end{cor}

\begin{proof}
By \cite{ASZ} Lemma 1.1 the categories
$\msf{M}_{q, \mrm{f}} / \msf{M}_{q - 1, \mrm{f}}$
and
$\msf{M}^{\circ}_{q, \mrm{f}} / \msf{M}^{\circ}_{q - 1, \mrm{f}}$
are noetherian, hence by Theorem \ref{a2.9} they are also artinian.
\end{proof}

At last here is:

\begin{proof}[Proof of Theorem \tup{\ref{a2.19}}]
We already verified axioms (i) and (ii).
The fact that $\opn{Cdim}_{R}$ is finitely partitive is immediate
from Corollary \ref{a2.10}, and this in turn easily implies
axiom (iii) -- cf.\ \cite{MR} Corollary 8.3.6.
\end{proof}

\begin{cor} \label{a2.12}
Every finitely generated $A$-module has a
$\opn{Cdim}_{R}$-critical composition series.
\end{cor}

\begin{proof}
By Theorem \ref{a2.19} and \cite{MR} Proposition 6.2.20.
\end{proof}

Let
\[ \begin{aligned}
d_{0} & := \inf \{ \opn{Cdim}_{R} M \mid M \in \cat{Mod} A,
M \neq 0 \} \\
d_{1} & := \opn{Cdim}_{R} A =
\sup \{ \opn{Cdim}_{R} M \mid M \in \cat{Mod} A \} .
\end{aligned} \]

\begin{cor} \label{a2.13}
Suppose
$R \in \msf{D}(\cat{Mod}(A \otimes B^{\circ}))$
is an Auslander dualizing complex. Then
\[ \opn{Kdim} M \leq \opn{Cdim}_{R} M - d_{0} \]
for all finitely generated $A$-modules $M$.
In particular if $\opn{Cdim}_{R} M = d_{0}$ then $M$ is artinian.
\end{cor}

\begin{proof}
By induction, starting with $q  = d_{0}$,
Theorem \ref{a2.6} and Corollary \ref{a2.10} show that
$\opn{Kdim} M \leq q - d_{0}$ for all $M \in \msf{M}_{q, \mrm{f}}$.
\end{proof}

Here is a generalization of \cite{Bj} Theorem 1.14.

\begin{thm} \tup{(Gabber's Maximality Principle)} \label{a3.1}
Let $A$ and $B$ be $k$-algebras, with $A$ left noetherian and $B$
right noetherian, and let
$R \in \msf{D}(\cat{Mod}(A \otimes B^{\circ}))$
be an Auslander dualizing complex.
Suppose $N$ is a $\opn{Cdim}_{R}$-pure $A$-module with
$\opn{Cdim}_{R} N = n$, and $M$ is a finitely generated submodule.
Then there is a unique maximal module $\tilde{M}$
such that $M \subset \tilde{M} \subset N$, $\tilde{M}$ is finitely
generated, and
$\opn{Cdim}_{R} \tilde{M} / M \leq n - 2$.
\end{thm}

\begin{proof}
Note that we do not assume $N$ is finitely generated.
The uniqueness is clear because $\opn{Cdim}_{R}$ is an exact
dimension function. So it remains to show existence.
If $\tilde{M}$ is any finitely generated submodule of $N$ containing
$M$, such that $\opn{Cdim}_{R} \tilde{M} / M \leq n - 2$, then
$\opn{Ext}^{-n}_{A}(M, R) \cong \opn{Ext}^{-n}_{A}(\tilde{M}, R)$.
Hence, by Theorem \ref{a2.6}(3), the module $\tilde{M}$ embeds
functorially into the finitely generated $A$-module
$\opn{Ext}^{-n}_{B^\circ} (\opn{Ext}^{-n}_{A}(M, R), R)$.
This implies there is a maximal such $\tilde{M}$.
\end{proof}

The next two theorems generalize \cite{GL} Theorems 1.4 and 1.6
by eliminating the Gorenstein and Cohen-Macaulay conditions.
From here on we consider a single noetherian algebra $A$
(i.e.\ $A = B$).

\begin{dfn} \label{a2.22}
Let $\opn{dim}$ be an exact dimension function.
\begin{enumerate}
\item $\opn{dim}$ is called {\em symmetric} if
$\opn{dim}_{A} M = \opn{dim}_{A^{\circ}} M$
for every bimodule $M$ finitely generated on both sides.
\item $\opn{dim}$ is called {\em weakly symmetric} if
$\opn{dim}_{A} M = \opn{dim}_{A^{\circ}} M$
for every bimodule $M$ which is a subquotient of $A$.
\end{enumerate}
\end{dfn}

\begin{lem} \label{a3.2}
Let $\opn{dim}$ be a weakly symmetric exact dimension function.
\begin{enumerate}
\item Let $M$ be a finitely generated $\opn{dim}$-pure
$A$-module, and let $I := \opn{Ann}_{A}(M)$. If
$\opn{dim} M = \opn{dim} A / I$ then $A / I$ is
$\opn{dim}$-pure.
\item Let $I$ be an ideal of $A$ such that $A / I$ is
$\opn{dim}$-pure and let $\mfrak{q}$ be a prime ideal of $A$
that is minimal over $I$. Then
$\opn{dim} A / \mfrak{q} = \opn{dim} A / I$.
\end{enumerate}
\end{lem}

\begin{proof}
This is completely analogous to \cite{KL} 9.6 and 9.5.
\end{proof}

Recall from \cite{GL} that $\opn{Spec} A$ is said to have {\em normal
separation} provided that for any pair prime ideals
$\mfrak{p} \subsetneqq \mfrak{q}$,
the factor $\mfrak{q} / \mfrak{p}$ contains a nonzero normal
element of $A / \mfrak{p}$. Under the assumptions of the lemma we say
that {\em Tauvel's height formula} holds in $A$ provided
\[ \opn{height} \mfrak{p} + \opn{dim} A / \mfrak{p} =
\opn{dim} A \]
for all primes $\mfrak{p}$.

\begin{thm} \label{a3.3}
Suppose that $A$ is a noetherian $k$-algebra, $R$ is an Auslander
dualizing complex over $A$ and $\opn{Cdim}_{R}$ is weakly symmetric.
Let $\mfrak{p} \subsetneqq \mfrak{q}$ be prime ideals of $A$ with
$\opn{height} \mfrak{q} / \mfrak{p} = 1$. If there exists an element
$a \in \mfrak{q} - \mfrak{p}$ that is
normal modulo $\mfrak{p}$, then
$\opn{Cdim}_{R} A / \mfrak{p} = \opn{Cdim}_{R} A / \mfrak{q} + 1$.
\end{thm}

\begin{proof}
Use the proof of \cite{GL} Theorem 1.4, but replace $\opn{GKdim}$ by
$\opn{Cdim}_{R}$, and use Lemma \ref{a3.2} instead of \cite{Len}
Lemmas 2 and 3.
\end{proof}

\begin{thm} \label{a3.4}
Suppose that $A$ is a noetherian $k$-algebra, $R$ is an Auslander
dualizing complex over $A$ and $\opn{Cdim}_{R}$ is weakly symmetric.
If $\opn{Spec} A$ is normally separated, then $A$ is catenary.
If in addition $A$ is prime, then Tauvel's height formula holds.
\end{thm}

\begin{proof}
The proof of \cite{GL} Theorem 1.6 works here after we replace
$\opn{GKdim}$ by $\opn{Cdim}_{R}$.
\end{proof}

We call attention to Question \ref{a4.15} regarding the possible
symmetry of $\opn{Cdim}_{R}$.

The Macaulay property of \cite{SZ} is adapted in the following way,
to be used later in the paper.

\begin{dfn} \label{a2.14}
Suppose $R$ is an Auslander dualizing complex over $A$.
Let $\opn{dim}$ be an exact dimension function on finitely generated
$A$-modules. If there is some integer $c$ such that
\[ \opn{dim} M = \opn{Cdim}_{R} M + c \]
for all $M \in \cat{Mod}_{\mrm{f}} A$,
then we say $R$ is {\em Macaulay with respect to $\opn{dim}$}, or
that $R$ is {\em $\opn{dim}$-Macaulay}.
\end{dfn}

Note that if $R = A$, then ``$\opn{GKdim}$-Macaulay''
is equivalent to ``Cohen-Macaulay'' as it is used in \cite{Bj},
\cite{Lev}, \cite{ASZ}, and \cite{SZ}. This is because
$\opn{GKdim} M +j_{R}(M) = c = \opn{dim} A$
in this case.

\begin{exa} \label{a2.15}
If $A$ is a commutative affine $k$-algebra
and $R$ is any central dualizing complex over $A$, then $R$ is
Auslander and $\opn{GKdim}$-Macaulay. In this case we also have
$\opn{Kdim} M = \opn{GKdim} M$ for all finitely generated $A$-modules
$M$.
\end{exa}

\section{Rigid Dualizing Complexes}

In this section we consider dualizing complexes which satisfy
a special condition discovered by Van den Bergh \cite{VdB}.
Rigid dualizing complexes are unique and even functorial. Furthermore
if $R$ is an Auslander rigid dualizing complex then the canonical
dimension $\opn{Cdim}_{R}$ is particularly well behaved (as examples
indicate; see Question \ref{a4.15}). By default $A$ and $B$ denote
noetherian $k$-algebras.

First we shall need some more notation for bimodules.
Suppose $A$ and $B$ are $k$-algebras.
For an element $a \in A$ we denote by $a^{\circ} \in A^{\circ}$
the same element. Thus for $a_{1}, a_{2} \in A$,
$a_{1}^{\circ} \cdot a_{2}^{\circ} =
(a_{2} \cdot a_{1})^{\circ} \in A^{\circ}$.
With this notation if $M$ is a right $A$-module then the left
$A^{\circ}$ action is
$a^{\circ} \cdot m = m \cdot a$, $m \in M$.
The algebra $A^{\mrm{e}}$ has an involution
$A^{\mrm{e}} \ar (A^{\mrm{e}})^{\circ}$,
$a_{1} \otimes a_{2}^{\circ} \mapsto a_{2} \otimes a_{1}^{\circ}$
which allows us to regard every left $A^{\mrm{e}}$-module $M$
as a right $A^{\mrm{e}}$-module in a consistent way:
\[ (a_{1} \otimes a_{2}^{\circ}) \cdot m =
(a_{2} \otimes a_{1}^{\circ})^{\circ} \cdot m =
m \cdot (a_{2} \otimes a_{1}^{\circ}) = a_{1} \cdot m \cdot a_{2} . \]

Given an $(A \otimes B)$-module $M$ and a $(B \otimes A)$-module $N$
we define a mixed action of
$A^{\mrm{e}} \otimes B^{\mrm{e}}$
on the tensor product $M \otimes N$ as follows.
$A^{\mrm{e}}$ acts on $M \otimes N$  by the outside action
\[ (a_{1} \otimes a_{2}^{\circ}) \cdot (m \otimes n) :=
(a_{1} \cdot m) \otimes (n \cdot a_{2}) , \]
whereas $B^{\mrm{e}}$ acts on $M \otimes N$  by the inside action
\[ (b_{1} \otimes b_{2}^{\circ}) \cdot (m \otimes n) :=
(m \cdot b_{2}) \otimes (b_{1} \cdot n) . \]
By default we regard the outside action as a left action and the
inside action as a right action. If $A = B$ and $M = N$ then
the two actions by $A^{\mrm{e}}$ on $M \otimes M$ are interchanged
the involution
$m_{1} \otimes m_{2} \mapsto m_{2} \otimes m_{1}$. However for
the sake of definiteness in this case, given an $A^{\mrm{e}}$-module
$L$,
$\opn{Hom}_{A^{\mrm{e}}}(L, M \otimes M)$ shall refer to
homomorphisms $L \ar M \otimes M$ which are $A^{\mrm{e}}$-linear
with respect to the outside action.

\begin{dfn}[\cite{VdB} Definition 8.1] \label{a4.1}
Suppose $R$ is a dualizing complex over $A$. If there is an
isomorphism
\[ \phi : R \iso \mrm{R} \opn{Hom}_{A^{\mrm{e}}}(A, R \otimes R) \]
in $\msf{D}(\cat{Mod} A^{\mrm{e}})$,
we call $(R, \phi)$ a {\em rigid dualizing complex}.
\end{dfn}

It is obvious that if $R$ is rigid, then any shift $R[n]$, for
$n \neq 0$, is no longer rigid.
Van den Bergh proved that a rigid dualizing complex $(R, \phi)$
over $A$ is unique, up to an isomorphism in
$\msf{D}(\cat{Mod} A^{\mrm{e}})$ (see \cite{VdB} Proposition 8.2).
Below we extend this result by proving that rigid dualizing complexes
are functorial, in a suitable sense.

Let $A \ar B$ be a $k$-algebra homomorphism.
Given $M \in \msf{D}^{+}(\cat{Mod} B)$,
$N \in \msf{D}^{+}(\cat{Mod} A)$ and a morphism
$\psi : M \ar N$ in $\msf{D}(\cat{Mod} A)$, $\psi$ factors
naturally through
$\mrm{R} \opn{Hom}_{A}(B, N)$. This can be seen by replacing $N$ with
an injective resolution $I$ in $\msf{D}^{+}(\cat{Mod} A)$, and then
we can take $\psi$ to be a homomorphism of complexes.
The image of $\psi$ will then land inside
$\opn{Hom}_{A}(B, I)$.
The same fact is true for bimodules.

We say a $k$-algebra homomorphism $A \ar B$
if {\em finite} if $B$ is finitely generated as a left and as a
right $A$-module.

\begin{thm} \label{a4.2}
Let $A \ar B$ be a finite $k$-algebra homomorphism.
Suppose $A$ and $B$ have rigid dualizing complexes
$(R_{A}, \phi_{A})$ and $(R_{B}, \phi_{B})$ respectively.
Then there is {\em at most one} morphism
$\psi : R_{B} \ar R_{A}$
in $\msf{D}(\cat{Mod} A^{\mrm{e}})$ satisfying conditions \tup{(i)}
and \tup{(ii)} below:
\begin{enumerate}
\rmitem{i} $\psi$ induces an isomorphism
\[ R_{B} \cong \mrm{R} \opn{Hom}_{A}(B, R_{A}) \cong
\mrm{R} \opn{Hom}_{A^{\circ}}(B, R_{A}) \]
in $\msf{D}(\cat{Mod} A^{\mrm{e}})$.
\rmitem{ii} The diagram
\[  \begin{CD}
R_{B} @> \phi_{B} >>
\mrm{R} \opn{Hom}_{B^{\mrm{e}}}
(B, R_{B} \otimes R_{B}) \\
@V \psi VV @V \psi \otimes \psi VV \\
R_{A} @> \phi_{A} >>
\mrm{R} \opn{Hom}_{A^{\mrm{e}}}
(A, R_{A} \otimes R_{A})
\end{CD} \]
in $\msf{D}(\cat{Mod} A^{\mrm{e}})$ is commutative.
\end{enumerate}
\end{thm}

The theorem is proved after this lemma.
Given a $k$-algebra homomorphism $A \ar B$, denote by
$\mrm{Z}_{B}(A) \subset B$ the centralizer of $A$.

\begin{lem} \label{a4.3}
Let $A \ar B$ be a finite $k$-algebra homomorphism.
Suppose $A$ and $B$ have dualizing complexes
$R_{A}$ and $R_{B}$ respectively, and
$\psi : R_{A} \ar R_{B}$
is a morphism in $\msf{D}(\cat{Mod} A^{\mrm{e}})$ satisfying
condition \tup{(i)} of the theorem. Then
$\opn{Hom}_{\msf{D}(\cat{Mod} A^{\mrm{e}})}(R_{B}, R_{A})$
is a free left and right $\mrm{Z}_{B}(A)$-module with basis $\psi$.
\end{lem}

\begin{proof}
Denote by
$D := \mrm{R} \opn{Hom}_{A}(-, R_{A})$
and
$D^{\circ} := \mrm{R} \opn{Hom}_{A^{\circ}}(-, R_{A})$
the dualizing functors. By assumption $R_{B} \cong D^{\circ} B$.
Applying the functor $D$ we get
isomorphisms of left $\mrm{Z}_{B}(A)$-modules
\[ \begin{aligned}
\opn{Hom}_{\msf{D}(\cat{Mod} A^{\mrm{e}})}(R_{B}, R_{A}) & \cong
\opn{Hom}_{\msf{D}(\cat{Mod} A^{\mrm{e}})}(D^{\circ} B, D^{\circ} A)
\\ & \cong
\opn{Hom}_{\msf{D}(\cat{Mod} A^{\mrm{e}})}(A, B) \\
& \cong \opn{Hom}_{A^{\mrm{e}}}(A, B) \cong \mrm{Z}_{B}(A) .
\end{aligned} \]
Likewise for the right action.
\end{proof}

\begin{proof}[Proof of the theorem]
Assume $\psi'$ is another such isomorphism. According to the
lemma above,
\[ \psi' = (b_{1} \otimes 1) \psi = (1 \otimes b_{2}^{\circ}) \psi \]
for suitable $b_{i} \in \mrm{Z}_{B}(A)^{\times}$. So
\[ \psi' \otimes \psi' = (b_{1} \otimes b_{2}^{\circ})
(\psi \otimes \psi) . \]
Now the diagram in condition (ii) consists of morphisms in
$\msf{D}(\cat{Mod} A^{\mrm{e}})$.
Since multiplications by $b_{1}$ and $b_{2}^{\circ}$ are
$A^{\mrm{e}}$-linear, we see that
\[ (b_{1} \otimes 1) \psi = (b_{1} \otimes b_{2}^{\circ}) \psi
\in \opn{Hom}_{\msf{D}(\cat{Mod} A^{\mrm{e}})}(R_{B}, R_{A}) . \]
Hence dividing by the unit $b_{1} \otimes 1$,
we see that $b_{2} = 1$.
\end{proof}

\begin{cor} \label{a4.6}
A rigid dualizing complex $(R, \phi)$ over $A$ is unique up to
a unique isomorphism.
\end{cor}

\begin{proof}
Suppose $(R', \phi')$ is another rigid dualizing complex.
We claim that there is an isomorphism $\psi : R \iso R'$ in
$\msf{D}(\cat{Mod} A^{\mrm{e}})$ which satisfies condition (ii)
of the theorem. By Theorem \ref{a4.2} such $\psi$ is unique.

To produce $\psi$, choose any isomorphism $\psi' : R \iso R'$,
which we know exists by \cite{VdB} Proposition 8.2.
Then by Lemma \ref{a4.3} there are $a_{i} \in \mrm{Z}(A)^{\times}$
such that
\[ (a_{1} \otimes 1)^{-1} \psi' =
(1 \otimes a_{2}^{\circ})^{-1} \psi' =
{\phi'}^{-1} (\psi' \otimes \psi') \phi . \]
The isomorphism
\[ \psi := (a_{1} \otimes 1) \psi' =
(1 \otimes a_{2}^{\circ}) \psi' \]
will satisfy condition (ii).
\end{proof}

Thus we may speak of {\em the} rigid dualizing complex of $A$ (if it
exists).

Lemma \ref{a4.3} can be sharpened when $A = B$.

\begin{prop} \label{a4.4}
Let $(R, \phi)$ be a rigid dualizing complex over $A$.
Then the two $k$-algebra homomorphisms
\[ \lambda, \rho : \mrm{Z}(A) \ar
\opn{End}_{\msf{D}(\cat{Mod} A^{\mrm{e}})}(R) \]
from the center of $A$, namely left and right multiplication,
are both isomorphisms, and are equal.
\end{prop}

\begin{proof}
By Lemma \ref{a4.3} with $A = B$ and $\psi = 1$, we see that
$\lambda$ and $\rho$ are isomorphisms. Take $a \in \mrm{Z}(A)$, and
let $a' := \rho^{-1} \lambda(a) \in \mrm{Z}(A)$.
Using the definition of the mixed action on $R \otimes R$ and
the rigidification isomorphism $\phi$, the commutation (or
conjugation) of $a$ across $R$ is transferred to commutation of $a'$
across $A$. Since $a'$ does commute with $A$ it follows that
$\lambda(a) = \rho(a)$ (and so in fact $a' = a$).
\end{proof}

We will often omit reference to the rigidifying isomorphism $\phi$.

\begin{cor} \label{a4.5}
If $R$ is a rigid dualizing complex over $A$ then for any
$q$ the cohomology bimodule $\mrm{H}^{q} R$ is central over
$\mrm{Z}(A)$.
\end{cor}

\begin{dfn} \label{a4.7}
Let $A \ar B$ be a finite homomorphism of $k$-algebras.
Assume the rigid dualizing complexes
$(R_{A}, \phi_{A})$ and $(R_{B}, \phi_{B})$ exist.
If there is a morphism $\psi$ satisfying the conditions of Theorem
\ref{a4.2} then we call it the {\em trace morphism}
and denote it by $\opn{Tr}_{B / A}$.
\end{dfn}

The next corollary is obvious.

\begin{cor} \label{a4.8}
Let $A \ar B$ and $B \ar C$ be finite $k$-algebra homomorphisms.
Assume the rigid dualizing complexes
$(R_{A}, \phi_{A})$, $(R_{B}, \phi_{B})$ and $(R_{C}, \phi_{C})$
and the trace morphisms $\opn{Tr}_{B / A}$ and $\opn{Tr}_{C / B}$
exist. Then
$\opn{Tr}_{C / A}$ exists too, and
\[ \opn{Tr}_{C / A} = \opn{Tr}_{B / A} \opn{Tr}_{C / B} . \]
\end{cor}

The existence of the trace morphism allows to transfer good
properties of $R_{A}$ to $R_{B}$.

\begin{prop} \label{a4.9}
Let $A \ar B$ be a finite homomorphism of $k$-algebras, and assume
the rigid dualizing complexes $R_{A}$ and $R_{B}$ and the trace
morphism $\opn{Tr}_{B / A}$ exist.
\begin{enumerate}
\item Let $C$ be any $k$-algebra. Then for
$M \in \msf{D}(\cat{Mod} (B \otimes C^{\circ}))$
there is a functorial isomorphism
\[ \mrm{R} \opn{Hom}_{B}(M, R_{B}) \cong
\mrm{R} \opn{Hom}_{A}(M, R_{A}) \]
in $\msf{D}(\cat{Mod} (C \otimes A^{\circ}))$.
\item If $R_{A}$ is Auslander then so is $R_{B}$.
\item For any $B$-module $M$,
$\opn{Cdim}_{R_{A}; A} M = \opn{Cdim}_{R_{B}; B} M$.
\item If $R_{A}$ is $\opn{GKdim}$-Macaulay, then
so is $R_{B}$.
\item Suppose $A \ar B$ is surjective. If $R_{A}$ is
$\opn{Kdim}$-Macaulay then so is $R_{B}$.
\end{enumerate}
\end{prop}

\begin{proof}
1. We can assume $R_{A}$ and $R_{B}$ are complexes of
injectives over $A^{\mrm{e}}$ and $B^{\mrm{e}}$ respectively,
and $\opn{Tr}_{B / A}$ is a homomorphism of complexes. Then we get a
functorial morphism
$\mrm{R} \opn{Hom}_{B}(M, R_{B}) \ar \mrm{R} \opn{Hom}_{A}(M, R_{A})$
in $\msf{D}(\cat{Mod} (C \otimes A^{\circ}))$.
To prove it's an isomorphism we can forget the $C$-module
structure. Because the two functors are way-out in both directions
(see \cite{RD} Section I.7), and they send direct sums to direct
products, it suffices to check for an isomorphism when $M = B$.
But that's given.\\
2, 3. Let $M$ be a finitely generated $B$-module and
$N \subset \opn{Ext}^{q}_{B}(M, R_{B})$ a $B^{\circ}$-submodule.
Then by part 1,
$N \subset \opn{Ext}^{q}_{A}(M, R_{A})$ as $A^{\circ}$-modules,
and for every $p$,
$\opn{Ext}^{p}_{B^{\circ}}(N, R_{B}) \cong
\opn{Ext}^{p}_{A^{\circ}}(N, R_{A})$
as $A$-modules. This proves the Auslander condition for $B$
and the dimension equality for finitely generated $B$-modules.\\
4. Follows from part 3 and the fact
$\opn{GKdim}_{A} M = \opn{GKdim}_{B} M$.\\
5. Similar to 4.
\end{proof}

\begin{dfn} \label{a4.11}
Suppose $A$ has an Auslander rigid dualizing complex $R$. Then we
denote the canonical dimension $\opn{Cdim}_{R}$ by $\opn{Cdim}$.
\end{dfn}

\begin{exa} \label{a4.10}
Suppose $A$ is an affine $k$-algebra and finite over its center.
Then we can find a smooth integral commutative $k$-algebra $C$ (e.g.\
a polynomial algebra),
and a finite homomorphism $C \ar \mrm{Z}(A)$. Say $\opn{Kdim} C = n$.
Since $\Omega^{n}_{C / k}[n]$ is a rigid dualizing complex over
$C$, it follows from \cite{Ye3} Proposition 5.9 that
\[ R := \mrm{R} \opn{Hom}_{C}(A, \Omega^{n}_{C / k}[n]) \]
is a rigid dualizing complex over $A$.
\end{exa}

\begin{exa}
Let $A$ be the algebra
$\left( \begin{smallmatrix} k & V \\ 0 & k \end{smallmatrix}
\right)$
where $V$ is a finite rank $k$-module. The rigid dualizing complex
is $A^{*} := \opn{Hom}_{k}(A, k)$.
Now $A$ is hereditary, hence Gorenstein, so the bimodule $A$ is a
dualizing complex. When $V \neq 0$ the dualizing
complexes $A$ and $A^{*}$ are not isomorphic, so $A$ is not rigid
then.
\end{exa}

\begin{exa}
Let $t_{1}, t_{2}, \ldots$ be a countable sequence of commuting
indeterminates
and let $C := k(t_{1}, t_{2}, \ldots)$ be the field of rational
functions. We claim that as $k$-algebra, $C$ has no rigid dualizing
complex. Since any dualizing complex over $C$ has to be of the form
$R = C^{\sigma}[n]$ for an automorphism $\sigma$ and an integer $n$
(by \cite{Ye3} Corollary 4.6 and Propositions 3.4 and 3.5), it
suffices to prove that
$\opn{Ext}^{i}_{C^{\mrm{e}}}(C, C^{\mrm{e}}) = 0$
for all $i$. This follows from the next lemma with $n = 0$.
\end{exa}

\begin{lem} \label{a5.8}
Let $D_{n} := C^{\mrm{e}} / I$ where $I$ is the ideal generated by
the elements
$f_{j} := x_{j} \otimes 1 - 1 \otimes x_{j}$ for $j = 1, \ldots, n$.
Then
$\opn{Ext}^{i}_{C^{\mrm{e}}}(C, D_{n}) = 0$
for all $i, n \geq 0$.
\end{lem}

\begin{proof}
Assume on the contrary that
$\opn{Ext}^{i}_{C^{\mrm{e}}}(C, D_{n}) \neq 0$ for some $i$ and $n$.
Let $i_0$ be the smallest such $i$. Since $f_{n + 1}$ is nonzero in
the domain $D_{n}$ there is a short exact sequence
\[ 0 \ar D_{n} \xrightarrow{f_{n + 1}} D_{n} \ar D_{n+1} \ar 0  \]
of $C^{\mrm{e}}$-modules.
That induces an exact sequence
\[ 0 = \opn{Ext}^{i_{0} - 1}_{C^{\mrm{e}}}(C, D_{n+1}) \ar
\opn{Ext}^{i_{0}}_{C^{\mrm{e}}}(C, D_{n})
\xrightarrow{f_{n + 1}}
\opn{Ext}^{i_{0}}_{C^{\mrm{e}}}(C, D_{n}) . \]
But $f_{n + 1}$ annihilates the $C^{\mrm{e}}$-module $C$,
which implies
$\opn{Ext}^{i_{0}}_{C^{\mrm{e}}}(C, D_{n}) = 0$,
contradicting the choice of $i_{0}$.
\end{proof}

We end the section with a basic question.

\begin{que} \label{a4.15}
Let $R$ be a rigid dualizing complex and $M$ an $A$-bimodule
finitely generated on both sides. Is there a functorial isomorphism
$\mrm{R} \opn{Hom}_{A}(M, R) \cong
\mrm{R} \opn{Hom}_{A^{\circ}}(M, R)$?
Or, at least, is $\opn{Cdim}_{R; A} M = \opn{Cdim}_{R; A^{\circ}} M$?
\end{que}

For a partial answer turn to Section 6, where the presence of
auxiliary filtrations allows us to take advantage
of results in Sections 4-5 on graded algebras.

\section{Dualizing Complexes over Graded Algebras}

In this section we consider connected $\mbb{Z}$-graded $k$-algebras,
namely algebras $A = \bigoplus_{n \geq 0} A_{n}$ with $A_{0} \cong k$
and $A_{n}$ a finitely generated $k$-module.

For such an algebra $A$ let $\cat{GrMod} A$ be the category of
$\mbb{Z}$-graded left modules, with degree $0$ homomorphisms.
For $M, N \in \cat{GrMod} A$ we write $M(n)$ for the shifted module
with $M(n)_{i} = M_{n + i}$, and
\[ \opn{Hom}_{A}^{\mrm{gr}}(M, N) := \bigoplus_{n \in \mbb{Z}}
\opn{Hom}_{\cat{GrMod} A}^{\mrm{gr}}(M, N(n)) . \]
There is a forgetful functor $\cat{GrMod} A \ar \cat{Mod} A$.
Observe that
$\opn{Hom}_{A}^{\mrm{gr}}(M, N) \subset \opn{Hom}_{A}(M, N)$
with equality when $M$ is finitely generated.

$\cat{GrMod} A$ is an abelian category with direct and inverse limits,
enough injectives and enough projectives.
Let $\cat{D}(\cat{GrMod} A)$ be the derived category.
The derived functors
$\mrm{R} \opn{Hom}^{\mrm{gr}}_{A}(M, N)$ and
$M {\otimes}^{\mrm{L}}_{A} N$
are calculated just as in the ungraded case, see Section 1, but using
graded-projectives or graded-injectives.

We say $M \in \cat{GrMod} k$ is {\em locally finite} if
each $M_{n}$ is a finitely generated $k$-module.
Let $M^{*} := \opn{Hom}_{k}^{\mrm{gr}}(M, k)$.
Denote by $\cat{D}_{\mrm{lf}}( \cat{GrMod} A)$ the subcategory
of complexes with locally finite cohomologies.
Matlis duality says that $M \cong M^{**}$ for
$M \in \cat{D}_{\mrm{lf}}(\cat{GrMod} A)$.

We denote by $\mfrak{m}$ the augmentation ideal
$\bigoplus_{n > 0} A_{n}$ of $A$, and we write
$\Gamma_{\mfrak{m}} M$
for the $\mfrak{m}$-torsion submodule of a graded $A$-module
$M$. There is a derived functor $\mrm{R} \Gamma_{\mfrak{m}}$,
which is calculated by graded-injectives (see \cite{Ye1}).
The cohomology modules are
$\mrm{H}^{i} \mrm{R} \Gamma_{\mfrak{m}} M =
\lim\limits_{n \ar} \opn{Ext}^{i}_{A}(A / A_{\geq n}, M)$.
We write $\mfrak{m}^{\circ}$ for the augmentation ideal of
$A^{\circ}$.

The definitions and results of the previous sections can all be
translated to the graded category by adding the adjective ``graded''
where needed, like ``graded dualizing complex,'' ``graded Auslander
property'' etc. The proofs of the graded variants of these results are
identical to the ungraded ones, so there is no need to repeat them.
In the rest of the paper we shall refer to such a result by writing
something like ``according to the graded variant of Theorem ...".

\begin{rem} \label{a6.1}
Let $G$ be any finitely generated abelian group.
Fix an isomorphism $G \cong \mbb{Z}^{r} \times T$, where $T$ is a
finite group, and a basis $g_{1}, \ldots, g_{r}$ of $\mbb{Z}^{r}$.
Let $G_{+}$ be the semigroup generated by $0$ and the elements
$g_{i} + t$, $1 \leq i \leq r$, $t \in T$.
For $g, g' \in G$ we write $g \geq g'$ if $g - g' \in G_{+}$,
and this defines a partial order on $G$.
A $G$-graded $k$-algebra $A$ is called {\em connected} if
$A = \bigoplus_{g \in G_{+}} A_{g}$,
$A_{0} = k$ and each $A_{g}$ is finitely generated as a module
over $k$. The augmentation ideal of $A$ is
$\mfrak{m} := \bigoplus_{g > 0} A_{g}$.

Note that the group homomorphism $\phi : G \ar \mbb{Z}$
sending $g_{i} \mapsto 1$
makes $A$ into a connected $\mbb{Z}$-graded algebra, with
$A_{n} = \bigoplus_{\phi(g) = n} A_{g}$, $n \in \mbb{Z}$.

It is not hard to see that all results in this paper which are stated
for connected $\mbb{Z}$-graded algebras are also valid for
connected $G$-graded $k$-algebras, for any $G$ as above.
\end{rem}

Throughout this section $A$ and $B$ are connected graded noetherian
$k$-algebras.

\begin{dfn}[\cite{Ye1} Definition 4.1] \label{a6.3}
A {\em balanced dualizing complex} over $A$
is a graded dualizing complex $R$ such that\
\[ \mrm{R} \Gamma_{\mfrak{m}} R \cong
\mrm{R} \Gamma_{\mfrak{m}^{\circ}} R \cong A^{*}  \]
in $\cat{D}(\cat{GrMod} A^{\mrm{e}})$.
\end{dfn}

A balanced dualizing complex $R$ is unique up
to isomorphism in $\cat{D}(\cat{GrMod} A^{\mrm{e}})$,
and its endomorphisms are just elements of $k$.
Thus if we choose an isomorphism
$\phi : \mrm{R} \Gamma_{\mfrak{m}} R \iso A^{*}$
in $\cat{D}(\cat{GrMod} A^{\mrm{e}})$, the pair
$(R, \phi)$ is unique up to a unique isomorphism.

It had been known for some time (by \cite{Ye1}) that balanced
dualizing complexes exist for Artin-Schelter Gorenstein algebras,
twisted homogeneous coordinate algebras
and algebras finite over their centers.
Recently additional existence results became available, due to the
work of Van den Bergh. First recall the following definition
taken from \cite{AZ}.

\begin{dfn} \label{a6.4}
The condition $\chi$ holds for a noetherian connected graded
$k$-algebra $A$ if for every $M \in \cat{GrMod}_{\mrm{f}} A$
and integer $i$, $\opn{Ext}^{i}_{A}(k, M)$
is a finitely generated $k$-module.
\end{dfn}

In view of \cite{AZ} Proposition 3.1(3), this definition is
equivalent to \cite{AZ} Definition 3.2; and by \cite{AZ}
Proposition 3.11(2) it is equivalent to \cite{AZ} Definition 3.7.
The next lemma provides further characterization of the condition
$\chi$. Recall that a graded module $M$ is said to be {\em right
bounded} if $M_{n} = 0$ for $n \gg 0$.

\begin{lem} \label{newlem4}
Let $A$ be a noetherian connected graded $k$-algebra and
$M \in \cat{GrMod}_{\mrm{f}} A$. Then the following are
equivalent:
\begin{enumerate}
\rmitem{i} $\opn{Ext}^{i}_{A}(k, M)$ is a finitely generated
$k$-module for all $i$.
\rmitem{ii} $\mrm{H}^{i} \mrm{R} \Gamma_{\mfrak{m}} M$
is right bounded for all $i$.
\rmitem{iii}
$\mrm{H}^{i} \mrm{R} \Gamma_{\mfrak{m}} M$
is an artinian $A$-module for all $i$.
\end{enumerate}
\end{lem}

\begin{proof}
(i) $\Leftrightarrow$ (ii) is by \cite{AZ} Corollary 3.6(3).
(iii) $\Rightarrow$ (ii) is immediate, since the socle of
$\mrm{H}^{i} \mrm{R} \Gamma_{\mfrak{m}} M$
is a finitely generated $k$-module, hence bounded.
Finally assume (i), and let $I$ be a minimal graded-injective
resolution of $M$. From \cite{Ye1} Lemma 4.3 it follows that
$\Gamma_{\mfrak{m}} I^{i} \cong A^{*} \otimes
\opn{Ext}^{i}_{A}(k, M)$,
which is artinian. Hence
$\mrm{H}^{i} \mrm{R} \Gamma_{\mfrak{m}} M$ is artinian.
\end{proof}

In an earlier paper we proved the next theorem.

\begin{thm}[\tup{\cite{YZ1} Theorem 4.2}] \label{ourthm}
Let $A$ be a noetherian connected graded $k$-algebra. If $A$ has a
balanced dualizing complex then the condition $\chi$ holds for $A$
and $A^{\circ}$, and the functors $\Gamma_{\mfrak{m}}$ and
$\Gamma_{\mfrak{m}^{\circ}}$ have finite cohomological dimensions.
\end{thm}

The converse, which is quite harder, was proved by Van den Bergh.

\begin{thm}[\tup{\cite{VdB} Theorem 6.3}] \label{a6.5}
Let $A$ be a noetherian connected graded $k$-algebra. Assume
the condition $\chi$ holds for $A$ and $A^{\circ}$, and the
functors $\Gamma_{\mfrak{m}}$ and $\Gamma_{\mfrak{m}^{\circ}}$
have finite cohomological dimensions. Then
\[ R_{A} := (\mrm{R} \Gamma_{\mfrak{m}} A)^{*} \]
is a balanced dualizing complex.
\end{thm}

Let us summarize some other known results related to
the balanced dualizing complexes.

\begin{thm}[Local Duality] \label{a6.6}
Let $R$ be a balanced dualizing complex over a noetherian connected
graded $k$-algebra $A$. Then for any graded $k$-algebra $B$ and any
$M \in \cat{D}(\cat{GrMod} (A \otimes B^{\circ}))$
there is a functorial isomorphism
\[ \mrm{R} \opn{Hom}^{\mrm{gr}}_{A}(M, R) \cong
(\mrm{R} \Gamma_{\mfrak{m}} M)^{*} . \]
\end{thm}

This is proved by combining \cite{VdB} Theorems 5.1 and 6.3.
The theorem was first proved in \cite{Ye1}, but only for
$M \in \cat{D}^{\mrm{b}}_{\mrm{f}}(\cat{GrMod} A)$.

\begin{prop}[\cite{VdB} Proposition 8.2(2)] \label{a6.25}
A balanced dualizing complex $R$ is rigid in the graded sense.
\end{prop}

\begin{rem} \label{a6.2}
According to an exercise in \cite{VdB} (whose only proof we know
is quite involved), if $I$ is a graded-injective $A$-module, then
$I$ has injective dimension $\leq 1$ in $\cat{Mod} A$. An immediate
consequence of this fact is that a graded dualizing complex
$R$ over $A$ is also an ungraded dualizing complex.
The special case we need, namely that a balanced dualizing
complex $R$ is rigid in the ungraded sense, is proved by other
means in Corollary \ref{a8.6}.
\end{rem}

Here is another result from \cite{VdB}.
Let us write $\mfrak{m}_{A^{\mrm{e}}}$ for the augmentation ideal
of $A^{\mrm{e}}$, so
$\mfrak{m}_{A^{\mrm{e}}} = \mfrak{m} \otimes A^{\circ} +
A \otimes \mfrak{m}^{\circ}$.

\begin{thm}[\cite{VdB} Corollary 4.8] \label{a6.7}
Assume $A$ has a balanced dualizing complex $R$. Let
$M \in \msf{D}(\cat{GrMod} A^{\mrm{e}})$ have finitely generated
cohomology modules on both sides. Then there is a functorial
isomorphism
\[ \mrm{R} \Gamma_{\mfrak{m}} M \cong
\mrm{R} \Gamma_{\mfrak{m}_{A^{\mrm{e}}}} M \cong
\mrm{R} \Gamma_{\mfrak{m}^{\circ}} M . \]
\end{thm}

We obtain the following interesting result:

\begin{cor} \label{a6.20}
Let $R$ be a balanced dualizing complex over $A$.
Then there is a functorial isomorphism
\[ \mrm{R} \opn{Hom}_{A}(M, R) \cong
\mrm{R} \opn{Hom}_{A^{\circ}}(M, R)  \]
for $M \in \msf{D}(\cat{GrMod} A^{\mrm{e}})$ with finitely generated
cohomology modules on both sides.
\end{cor}

We shall write $\opn{Cdim}_{A}$ instead of $\opn{Cdim}_{R; A}$
when $R$ is the balanced dualizing complex, and when we are
working in $\cat{GrMod} A$. Since a balanced dualizing complex
is rigid in the ungraded sense (by Corollary \ref{a8.6}), this is
consistent with Definition \ref{a4.11}.

\begin{dfn} \label{a6.21}
If $A$ has a graded Auslander balanced dualizing complex $R$ we say
$A$ is {\em graded Auslander}. Furthermore if $\opn{dim}$ is an exact
dimension function on graded modules, and if $R$ is graded
$\opn{dim}$-Macaulay, then we say $A$ is {\em graded Auslander
$\opn{dim}$-Macaulay}.
\end{dfn}

According to \cite{Ye1} Theorem 3.9, any two graded
dualizing complexes $R_{1}, R_{2}$ satisfy
$R_{2} \cong R_{1} \otimes_{A} A^{\sigma}(m)[n]$,
for an automorphism $\sigma$ and an integers $n, m$.
It follows that $R_{1}$ is graded
Auslander iff $R_{2}$ is so. In particular, if $A$ is graded
Gorenstein, then $A$ is graded Auslander-Gorenstein in the
usual sense iff it is graded Auslander in the sense of
Definition \ref{a6.21}.

Taking cohomologies in the previous corollary we get:

\begin{cor} \label{a6.9}
Suppose $A$ is graded Auslander. Then $\opn{Cdim}$ is symmetric
on graded modules. That is to say, if $M$ is a graded $A$-bimodule,
finitely generated on both sides, then
$\opn{Cdim}_{A} M = \opn{Cdim}_{A^{\circ}} M.$
\end{cor}

If $A$ is graded Auslander we have a bound on Krull dimension
of graded modules:

\begin{thm} \label{a6.11}
Suppose $A$ is graded Auslander. Then
\[ \opn{Kdim} M \leq \opn{Cdim} M=
\opn{sup} \{ q \mid \mrm{H}^{q} \mrm{R} \Gamma_{\mfrak{m}} M \neq 0 \}
< \infty  \]
for all finitely generated graded $A$-modules $M$.
\end{thm}

\begin{proof}
By Theorem \ref{a6.6}, if $M \neq 0$, we have
$\opn{Cdim} M =
\opn{sup} \{ q \mid \mrm{H}^{q} \mrm{R} \Gamma_{\mfrak{m}} M
\neq 0 \}$.
Next for a finitely generated graded $A$-module $M$
the Krull dimension is the same when computed in $\cat{GrMod} A$
and in $\cat{Mod} A$. And by the graded variant of Corollary
\ref{a2.13} we get
$\opn{Kdim} M \leq \opn{Cdim}_{R} M$,
since clearly $d_{0} = 0$.
\end{proof}

\begin{lem} \label{a6.13}
Let $A \ar B$ be a finite homomorphism of noetherian connected
graded $k$-algebras, with augmentation ideals
$\mfrak{m}_{A}, \mfrak{m}_{B}$. Assume $A$ satisfies condition
$\chi$. Then there is a functorial isomorphism
\[ \mrm{R} \Gamma_{\mfrak{m}_{B}} M \cong
\mrm{R} \Gamma_{\mfrak{m}_{A}} M \]
for $M \in \cat{D}^{+}(\cat{GrMod} B)$.
\end{lem}

\begin{proof}
For any homomorphism $A \ar B$ of graded algebras there is a
functorial morphism
$\mrm{R} \Gamma_{\mfrak{m}_{B}} M \ar
\mrm{R} \Gamma_{\mfrak{m}_{A}} M$
in $\cat{D}(\cat{GrMod} A)$.
By \cite{AZ} Lemma 8.2, the extra assumptions guarantee that
$\mrm{H}^{p} \mrm{R} \Gamma_{\mfrak{m}_{B}} M \ar
\mrm{H}^{p} \mrm{R} \Gamma_{\mfrak{m}_{A}} M$
is bijective for all $p$.
\end{proof}

The following theorem is a generalization of \cite{Jo1} Theorem 3.3.

\begin{thm} \label{a6.14}
Let $A \ar B$ be a finite homomorphism of graded $k$-algebras
and let $R_{A}$ be a balanced dualizing complex over $A$. Then:
\begin{enumerate}
\item $B$ has a balanced dualizing complex $R_{B}$.
\item There is a morphism
$\opn{Tr}_{B / A} : R_{B} \ar R_{A}$
in $\cat{D}(\cat{GrMod} A^{\mrm{e}})$,
which satisfies conditions \tup{(i)} and \tup{(ii)}
of Theorem \tup{\ref{a4.2}}.
\end{enumerate}
\end{thm}

\begin{proof}
From \cite{VdB} Theorem 6.3 we know that $A$ and $A^{\circ}$
satisfy $\chi$, and $\Gamma_{\mfrak{m}_{A}}$ and
$\Gamma_{\mfrak{m}^{\circ}_{A}} =\Gamma_{\mfrak{m}_{A^{\circ}}}$
have finite cohomological
dimensions. So by Lemma \ref{a6.13} the same is true for $B$.
Thus $B$ has a balanced dualizing complex
$R_{B} \cong (\mrm{R} \Gamma_{\mfrak{m}_{B}} B)^{*}$.

The morphism $A \ar B$ in $\cat{D}(\cat{GrMod} A^{\mrm{e}})$
induces a morphism
$(\mrm{R} \Gamma_{\mfrak{m}_{A}} B)^{*} \ar
(\mrm{R} \Gamma_{\mfrak{m}_{A}} A)^{*}$,
also in $\cat{D}(\cat{GrMod} A^{\mrm{e}})$. But
$(\mrm{R} \Gamma_{\mfrak{m}_{A}} B)^{*} \cong
(\mrm{R} \Gamma_{\mfrak{m}_{B}} B)^{*}$
and we get
$\opn{Tr}_{B / A} : R_{B} \ar R_{A}$.
The isomorphism of functors
$\mrm{R} \Gamma_{\mfrak{m}_{A}} \cong
\mrm{R} \Gamma_{\mfrak{m}_{A}^{\circ}}$
of \cite{VdB} Corollary 4.8 shows that $\opn{Tr}_{B / A}$ is the
same when calculated on the right,
i.e.,\ using $\mrm{R} \Gamma_{\mfrak{m}_{A}^{\circ}}$.

Condition (i) of Theorem \ref{a4.2} is a consequence of local
duality. To verify condition (ii) we again view $A \ar B$
as a morphism in $\cat{D}(\cat{GrMod} A^{\mrm{e}})$.
By \cite{VdB} Theorems 4.7 and 5.1 we get a commutative diagram
\[ \begin{CD}
(\mrm{R} \Gamma_{\mfrak{m}_{B}} B)^{*} @> \cong >>
(\mrm{R} \Gamma_{\mfrak{m}_{B^{\mrm{e}}}} B)^{*} @> \cong >>
\mrm{R} \opn{Hom}_{B^{\mrm{e}}}(B,
(\mrm{R} \Gamma_{\mfrak{m}_{B^{\mrm{e}}}} B^{\mrm{e}})^{*}) \\
@VVV @VVV @VVV \\
(\mrm{R} \Gamma_{\mfrak{m}_{A}} A)^{*} @> \cong >>
(\mrm{R} \Gamma_{\mfrak{m}_{A^{\mrm{e}}}} A)^{*} @> \cong >>
\mrm{R} \opn{Hom}_{A^{\mrm{e}}}(A,
(\mrm{R} \Gamma_{\mfrak{m}_{A^{\mrm{e}}}} A^{\mrm{e}})^{*})
\end{CD} \]
Finally by \cite{VdB} Theorem 7.1
\[ (\mrm{R} \Gamma_{\mfrak{m}_{A^{\mrm{e}}}} A^{\mrm{e}})^{*}
\cong (\mrm{R} \Gamma_{\mfrak{m}_{A}} A)^{*} \otimes
(\mrm{R} \Gamma_{\mfrak{m}_{A}^{\circ}} A)^{*} , \]
and of course the same for $B$.
\end{proof}

Applying the graded variant of Proposition \ref{a4.9} we obtain
the following corollary.

\begin{cor} \label{a6.15}
Let $A$ and $B$ be as in Theorem \tup{\ref{a6.14}}.
\begin{enumerate}
\item There is a functorial isomorphism
\[ \mrm{R} \opn{Hom}^{\mrm{gr}}_{B}(M, R_{B}) \cong
\mrm{R} \opn{Hom}^{\mrm{gr}}_{A}(M, R_{A}) \]
for all $M \in \cat{D}(\cat{GrMod} B)$.
\item If $A$ is graded Auslander then so is $B$.
\item $\opn{Cdim}_{A} M = \opn{Cdim}_{B} M$ for
$M \in \cat{GrMod} B$.
\item If $A$ is graded Auslander $\opn{GKdim}$-Macaulay then so is
$B$.
\item Suppose $A \ar B$ is surjective. If $A$ is graded Auslander
$\opn{Kdim}$-Macaulay, then so is $B$.
\end{enumerate}
\end{cor}

The next three propositions show that the graded Auslander property
can be transferred from one algebra to a related algebra.

\begin{prop} \label{a6.16}
Suppose $A$ has a balanced dualizing complex. Let $\opn{dim}$ stand
for either $\opn{Kdim}$ or $\opn{GKdim}$.
\begin{enumerate}
\item Let $\mfrak{a}, \mfrak{b}$ be graded ideals. If the
quotient algebras $A / \mfrak{a}$ and $A / \mfrak{b}$ are graded
Auslander \tup{(}resp.\ and graded $\opn{dim}$-Macaulay\tup{)},
then so is $A / \mfrak{a} \mfrak{b}$.
\item Let $\mfrak{a}$ be a nilpotent graded ideal of $A$.
If $A / \mfrak{a}$ is graded
Auslander \tup{(}resp.\ and graded $\opn{dim}$-Macaulay\tup{)}
then so is $A$.
\item If for every minimal graded prime ideal $\mfrak{p}$
the quotient algebra $A / \mfrak{p}$ is graded
Auslander  \tup{(}resp.\ and graded $\opn{dim}$-Macaulay\tup{)},
then so is $A$.
\end{enumerate}
\end{prop}

\begin{proof}
1. As usual we write
$D := \mrm{R} \opn{Hom}^{\mrm{gr}}_{A}(- , R)$
where $R$ is the balanced dualizing complex.
We may assume $\mfrak{a} \mfrak{b} = 0$.
Given a finitely generated graded module $M$ consider the exact
sequence
$0 \ar \mfrak{b} M \ar M \ar M / \mfrak{b} M \ar 0$,
and note that $\mfrak{b} M$ is an $A / \mfrak{a}$-module.
For any $i$ there is an exact sequence
$\mrm{H}^{i} D (M / \mfrak{b} M) \ar
\mrm{H}^{i} D M \ar \mrm{H}^{i} D (\mfrak{b} M)$.
Since $A / \mfrak{a}$ and $A / \mfrak{b}$ have the graded Auslander
property, the subquotients of $\mrm{H}^{i} D (M / \mfrak{b} M)$
and $\mrm{H}^{i} D (\mfrak{b} M)$
have $\opn{Cdim}$ no more than $i$.
Observe that here we are using Corollary \ref{a6.15}.
Hence by the long exact sequence
of duality, submodules of $\mrm{H}^{i} D M$ have $\opn{Cdim}$ no
more than $i$. The assertion about the Macaulay property is clear.\\
2. Use part 1 and induction.\\
3. Let $\mfrak{p}_{1}, \ldots, \mfrak{p}_{m}$ be the minimal
prime ideals of $A$. Then $(\prod_{i} \mfrak{p}_{i})^{n} = 0$
for some $n$, and we can use parts 1 and 2.
\end{proof}

\begin{prop} \label{a6.17}
Let $A \ar B$ be a finite homomorphism of connected graded algebras,
and assume $B \cong A \oplus L$ as graded $A$-bimodules.
If $B$ has a balanced dualizing complex then so does $A$.
\end{prop}

\begin{proof}
Let $M$ be a finitely generated graded $A$-module. Then by Lemma
\ref{a6.13} we get
\[ \mrm{R} \Gamma_{\mfrak{m}_{B}} (B \otimes_{A} M) \cong
\mrm{R} \Gamma_{\mfrak{m}_{A}} (B \otimes_{A} M) \cong
\mrm{R} \Gamma_{\mfrak{m}_{A}} M \oplus
\mrm{R} \Gamma_{\mfrak{m}_{A}} (L \otimes_{A} M) . \]
By Theorem \ref{ourthm} and Lemma \ref{newlem4} applied to $B$,
we see that the graded $B$-module
$\mrm{H}^{i} \mrm{R} \Gamma_{\mfrak{m}_{B}} (B \otimes_{A} M)$
is right bounded, and vanishes for large $i$.
Hence the same is true for the $A$-module
$\mrm{H}^{i} \mrm{R} \Gamma_{\mfrak{m}_{A}} M$.
Now apply Theorem \ref{a6.5}.
\end{proof}

We do not know if the under the assumptions of the proposition above
the Auslander property can be transferred from $B$ to $A$. However as
shown to us by Van den Bergh this is true in a special case:

\begin{prop}\label{a6.18}
Let $G$ be a finite group of order prime to $\opn{char} k$, acting
on $B$ by graded $k$-algebra automorphisms, and let
$A := B^{G}$ be the fixed ring.
If $B$ is graded Auslander then so is $A$.
\end{prop}

\begin{proof}
Given a finitely generated graded $A$-module $M$ and
a graded $A^{\circ}$-submodule
$N \subset \opn{Ext}^{q}_{A}(M, R_{A})$
we want to prove that
$\opn{Ext}^{p}_{A^{\circ}}(N, R_{A}) = 0$
for all $p < q$.
Let
$L := \{ b \in B \mid \sum_{g \in G} g(b) = 0 \}$,
so $B = A \oplus L$.
We have isomorphisms of graded $A^{\circ}$-modules
\begin{multline*}
\opn{Ext}^{q}_{A}(M, R_{A}) \oplus
\opn{Ext}^{q}_{A}(L \otimes_{A} M, R_{A}) \cong \\
\opn{Ext}^{q}_{A}(B \otimes_{A} M, R_{A}) \cong
\opn{Ext}^{q}_{B}(B \otimes_{A} M, R_{B})
\end{multline*}
that respect the $G$-action. Note that $G$ acts trivially on
$\opn{Ext}^{q}_{A}(M, R_{A})$.
Consider the graded $B^{\circ}$-module
$N \cdot B \subset \opn{Ext}^{q}_{B}(B \otimes_{A} M, R_{B})$.
Clearly
$N \cdot B = N + N \cdot L$.
But if
\[ n = \sum_{i} n_{i} l_{i} \in N \cap (N \cdot L) \]
with $n_{i} \in N$ and $l_{i} \in L$, then
\[ n = |G|^{-1} \sum_{g \in G} g \left(\sum_{i} n_{i} l_{i} \right)
= \sum_{i} n_{i} \cdot |G|^{-1} \sum_{g \in G} g(l_{i}) = 0 . \]
We conclude that
$N \cdot B = N \oplus N \cdot L$
as graded $A^{\circ}$-modules. Therefore
\begin{multline*}
\opn{Ext}^{p}_{A^{\circ}}(N, R_{A}) \oplus
\opn{Ext}^{p}_{A^{\circ}}(N \cdot L, R_{A}) \cong \\
\opn{Ext}^{p}_{A^{\circ}}(N \cdot B, R_{A}) \cong
\opn{Ext}^{p}_{B^{\circ}}(N \cdot B, R_{B}) = 0 .
\end{multline*}
\end{proof}

\section{Graded Algebras with Some Commutativity Hypothesis}

In this section we continue the discussion of balanced dualizing
complexes over connected graded noetherian $k$-algebras, but now we
look at algebras which have some commutativity hypothesis,
like PI, FBN, or enough normal elements. The main result here is:

\begin{thm} \label{a7.1}
Let $A$ be a noetherian connected graded $k$-algebra.
Suppose $t \in A$ is a homogeneous normal element of
positive degree, and let $B  := A / (t)$.
\begin{enumerate}
\item If $B$ has a balanced dualizing complex, then so does $A$.
\item If in addition $B$ is graded Auslander, then so is $A$.
\item If in addition $B$ is graded $\opn{Kdim}$-\tup{(}resp.\
$\opn{GKdim}$\tup{)}-Macaulay, then so is  $A$.
\end{enumerate}
\end{thm}

The ``classical'' case of this theorem, i.e.\ when $B$ is Gorenstein
and $t$ is a regular element (i.e.\ a non-zero-divisor), is
\cite{Lev} Theorem 3.6. Part 1 of the theorem is a trivial consequence
of Theorems \ref{ourthm} and \ref{a6.5}, and \cite{AZ} Theorem 8.8.

The proof of parts 2 and 3 appears after
a series of lemmas. In these lemmas we assume that $A$ has a balanced
dualizing complex (by part 1) and $B$ is graded Auslander.
The modules $M, N, \ldots$ will be finitely generated graded by
default. By Proposition \ref{a6.16} we can assume $A$ is prime, hence
$t$ is a regular element. The same proposition tells us that
\begin{equation} \label{a7.20}
A / (t^{n}) \text{ is graded Auslander for all } n \geq 1 .
\end{equation}
We denote by $D$ the duality functor
$\mrm{R} \opn{Hom}_{A}(-, R)$, where $R$ is the balanced dualizing
complex of $A$. Recall that according to Theorem \ref{a6.6},
$\mrm{H}^{-i} D M \cong (\mrm{R}^{i} \Gamma_{\mfrak{m}} M)^{*}$.
By definition
$\opn{Cdim} M = \sup \{ i \mid \mrm{H}^{-i} D M \neq 0 \}$,
so trivially
\begin{equation} \label{a7.21}
\opn{Cdim} M \leq \max \{ \opn{Cdim} M', \opn{Cdim} M / M' \}
\end{equation}
for all $M' \subset M$.

\begin{lem} \label{a7.2}
If $M$ is $t$-torsion-free then
$\opn{Cdim} M = \opn{Cdim} M / tM + 1$.
If $d= \opn{Cdim} M$ then $\mrm{H}^{-d} D M$ is $t$-torsion-free.
\end{lem}

\begin{proof}
Let $\sigma$ be the automorphism of $A$ such that
$t \cdot a = \sigma(a) \cdot t$, and let ${}^{\sigma} M$ be the
corresponding twisted module.
Then we have an exact sequence
\[  0 \ar {}^{\sigma^{-1}} M(-l) \xrightarrow{t \cdot} M
\ar M / t M \ar 0 \]
where $l$ is the degree of $t$. It is easy to see that
$\mrm{H}^{i} D ({}^{\sigma^{-1}} M(-l)) \cong
(\mrm{H}^{i} D M)^{\sigma^{-1}}(l)$,
so there is a long exact sequence
\begin{equation} \label{a7.3}
\mrm{H}^{i} D (M / t M) \ar \mrm{H}^{i} D M
\xrightarrow{\cdot t} (\mrm{H}^{i} D M)^{\sigma^{-1}}(l) \ar
\mrm{H}^{i + 1} D (M / t M) .
\end{equation}
If $\mrm{H}^{i + 1} D (M / t M) = 0$ then by the graded Nakayama
Lemma we get $\mrm{H}^{i} D M = 0$.
Therefore
\[ \opn{Cdim} M / t M \leq \opn{Cdim} M \leq \opn{Cdim} M / t M + 1
. \]

Now let $d = \opn{Cdim} M$. We need to show that
$\mrm{H}^{-d} D (M / t M) = 0$. If not, then
$\opn{Cdim} M / t M = d$, and hence also
$\opn{Cdim} M / t^{n} M = \opn{Cdim} t^{n - 1} M / t^{n} M = d$
for all $n \geq 1$.
According to Proposition \ref{a6.16} the algebra
$A / (t)^{n}$ has the graded Auslander property.
This implies that
$$\opn{Cdim} \mrm{H}^{-d} D (t^{n - 1} M / t^{n} M) = d
\quad {\text{and}}\quad
\opn{Cdim} \mrm{H}^{-d + 1} D (M / t^{n - 1} M) \leq d - 1.$$
Looking at the exact sequence
\begin{multline*}
0 \ar \mrm{H}^{-d} D (M / t^{n - 1} M) \ar
\mrm{H}^{-d} D (M / t^{n} M) \\
\ar \mrm{H}^{-d} D (t^{n - 1} M / t^{n} M) \xrightarrow{\phi}
\mrm{H}^{-d + 1} D (M / t^{n - 1} M)
\end{multline*}
we see that $\phi$ cannot be an injection.
Therefore
\[ \mrm{H}^{-d} D (M / t^{n - 1} M) \subsetneqq
\mrm{H}^{-d} D (M / t^{n} M) \subset
\mrm{H}^{-d} D (M) . \]
But this is true for all $n \geq 1$, contradicting with the noetherian
property of $\mrm{H}^{-d} D (M)$.
The upshot is that $\mrm{H}^{d} D (M / t M) = 0$.
Taking $i = -d$ in (\ref{a7.3}) we conclude that
$\mrm{H}^{d} D M$ is $t$-torsion-free.
\end{proof}

\begin{lem} \label{a7.4}
Let $N$ be the $t$-torsion submodule of $M$. Then
\[ \opn{Cdim} M = \max \{ \opn{Cdim} N, \opn{Cdim} M / N \} . \]
\end{lem}

\begin{proof}
Let $d$ be the right hand side.
Trivially $\opn{Cdim} M \leq d$ holds.
Assume $\opn{Cdim} M < d$. Dualizing the exact sequence
$0 \ar N \ar M \ar M / N \ar 0$
we obtain a long exact sequence
\[ \cdots \ar \mrm{H}^{i} D (M / N)  \ar \mrm{H}^{i} D M
\ar \mrm{H}^{i} D N \ar \mrm{H}^{i + 1} D (M / N)
\ar \mrm{H}^{i + 1} D M \ar \cdots \]
and taking $i = -d - 1$ we see that
$\mrm{H}^{-d} D (M / N) = 0$,  so $\opn{Cdim} M / N < d$.
Hence $\opn{Cdim} N = d$.
Also taking $i = -d$ we see that
$\mrm{H}^{-d} D N \subset \mrm{H}^{-d + 1} D (M / N)$.
Now $\mrm{H}^{-d} D N$ is $t$-torsion, yet by Lemma \ref{a7.2},
$\mrm{H}^{-d + 1} D (M / N)$ is $t$-torsion-free.
We conclude that $\mrm{H}^{-d} D N = 0$, which is a contradiction.
\end{proof}

\begin{lem} \label{a7.5}
\[ \opn{Cdim} M = \max \{ \opn{Cdim} L, \opn{Cdim} M / L \}  \]
for every $L \subset M$.
\end{lem}

\begin{proof}
By (\ref{a7.21}) it remains to prove ``$\geq$''.
Let $N$ be the $t$-torsion submodule of $M$. By (\ref{a7.20}) we have
$\opn{Cdim} N \geq \opn{Cdim} N \cap L$. On the other hand by Lemma
\ref{a7.4},
$\opn{Cdim} M = \max \{ \opn{Cdim} N, \opn{Cdim} M / N \}$
and
$\opn{Cdim} L = \max \{ \opn{Cdim} N \cap L, \opn{Cdim}$ \linebreak
$L / N \cap L \}$.
Hence it suffices to prove that
$\opn{Cdim} M / N = \max \{ \opn{Cdim} L / N \cap L, \opn{Cdim}$
\linebreak $M / L \}$.
So we may assume $M = M / N$ is $t$-torsion-free.

If $P := M / L$ is also $t$-torsion-free then we get a short exact
sequence
\[ 0 \ar L / t L \ar M / t M \ar P / t P \ar 0 . \]
Hence the assertion follows from Lemma \ref{a7.2}.

In general let $L' \supset L$ be such that $L' / L$ is the
$t$-torsion submodule of $M / L$. So $t^{n} L' \subset L$ for
some $n$. By Lemma \ref{a7.2},
\[ \opn{Cdim} L' / L \leq \opn{Cdim} L' / t^{n} L' =
\opn{Cdim} L' - 1 \]
and hence, from the long exact sequence of duality, we get
$\opn{Cdim} L = \opn{Cdim} L'$.
Applying the previous paragraph and Lemma \ref{a7.4} we have
\[ \begin{aligned}
\opn{Cdim} M &=\max \{ \opn{Cdim} L', \opn{Cdim} M / L' \} \\
& = \max \{ \opn{Cdim} L, \opn{Cdim} M / L' \} \\
& = \max \{ \opn{Cdim} L, \opn{Cdim} L' / L, \opn{Cdim} M / L' \} \\
& = \max \{ \opn{Cdim} L, \opn{Cdim} M / L \}
\end{aligned} \]
and we finish the proof.
\end{proof}

\begin{lem} \label{a7.6}
Suppose $t^{n}$ kills the $t$-torsion submodule of $M$. Then
\[ \opn{Cdim} M \leq \opn{Cdim} M / t^{n} M + 1 . \]
\end{lem}

\begin{proof}
Let $N$ be the $t$-torsion submodule of $M$ and $P := M / N$.
Then $t^{n} M \cong t^{n} P$ is an $t$-torsion-free submodule of $M$
and
$P / t^{n} P \cong M / (N \oplus t^{n} M)$.
By Lemmas \ref{a7.2} and \ref{a7.4},
\[ \opn{Cdim} M = \max \{ \opn{Cdim} N, \opn{Cdim} P \} =
\max \{ \opn{Cdim} N, \opn{Cdim} P / t^{n} P + 1 \} . \]
On the other hand by Lemma \ref{a7.5},
\[ \opn{Cdim} M / t^{n} M =
\max \{ \opn{Cdim}  N, \opn{Cdim} P / t^{n} P \} . \]
It remains to combine these equalities.
\end{proof}

\begin{proof}[Proof of Theorem \tup{\ref{a7.1}}]
As mentioned above part 1 is a consequence
of Theorems \ref{ourthm} and \ref{a6.5}, and \cite{AZ} Theorem 8.8.

\medskip \noindent 2.
Recall that the Auslander condition says that if
$N \subset \mrm{H}^{-i} D M$ then $\opn{Cdim} N \leq i$.
We can assume $A$ is prime and $t$ is regular.
By Lemma \ref{a7.5} it suffices to prove that
$\opn{Cdim} \mrm{H}^{-i} D M \leq i$ for all $i$. According to
(\ref{a7.20}) the inequality holds for $t$-torsion modules, so using
the long exact sequence of duality we may assume that $M$ is
$t$-torsion-free.

Choose $n$ such that \ $t^{n}$ kills the $t$-torsion submodule of
$\mrm{H}^{-i} D M$. The short exact sequence
\[ 0 \ar M \xrightarrow{t^{n} \cdot}
{}^{\sigma^{n}} M(n l)  \ar {}^{\sigma^{n}}(M / t^{n} M)(n l) \ar 0 \]
gives rise to a long exact sequence
\[ \cdots \ar (\mrm{H}^{-i} D M)^{\sigma^{n}}(-n l)
\xrightarrow{t^{n} \cdot} \mrm{H}^{-i} D M \ar
P = (\mrm{H}^{-i + 1} D (M / t^{n} M))^{\sigma^{n}}(-n l)
\ar \cdots . \]
Since $M / t^{n} M$ is $t$-torsion, $P$ is also $t$-torsion,
and every submodule of $P$ has $\opn{Cdim} \leq i - 1$.
So according to Lemma \ref{a7.6},
$\opn{Cdim} \mrm{H}^{-i} D M \leq (i - 1) + 1=i$.

\medskip \noindent 3.
Assume $B$ is graded $\opn{GKdim}$-Macaulay.
Since $A / (t^{n})$ is graded $\opn{GKdim}$-Macaulay for $n \geq 1$
(by Proposition \ref{a6.16}), it suffices to show that
$\opn{Cdim} M =\opn{GKdim} M$ for $t$-torsion-free modules. We know
that in this case
$\opn{GKdim} M = \opn{GKdim} M / t M + 1$, and by Lemma \ref{a7.2}
this is true also for $\opn{Cdim} M$.

Finally assume $B$ is graded $\opn{Kdim}$-Macaulay.
Again we need only consider $M$ which is $t$-torsion-free.
It is clear that
$\opn{Kdim} M \geq \opn{Kdim} M / t M + 1$.
Hence it follows from Lemmas \ref{a7.2} and \ref{a7.5} that
$\opn{Kdim} M \geq \opn{Cdim} M$. But by Theorem \ref{a6.11},
$\opn{Kdim} M \leq \opn{Cdim} M$.
\end{proof}

\begin{exa} \label{a7.7}
Let us consider a simple case of Theorem 7.1.
Let $A = B[t]$ where $t$ is a central variable of degree 1.
Let $R_{B}$ be the balanced dualizing complex of $B$. Then the
balanced dualizing complex of $A$ is nothing but
\[ R_{A} = R_{B} \otimes_{k[t]} \Omega^{1}_{k[t] / k}[1] \cong
R_{B}[t](-1)[1] \]
where
$\Omega^{1}_{k[t] / k} = k[t] \cdot \mrm{d} t \cong k[t](-1)$
as graded $k[t]$-modules. This follows from \cite{VdB} Theorem 7.1,
which states that if $B, C$ are noetherian connected graded
$k$-algebras with balanced dualizing complexes $R_{B}$ and $R_{C}$
respectively, and if $B \otimes C$ is noetherian, then
$R_{B} \otimes R_{C}$ is a balanced dualizing complex over
$B \otimes C$. By Theorem \ref{a7.1}, $A$ is graded Auslander and
$\opn{Kdim}$-(resp.\ $\opn{GKdim}$-)Macaulay if and only if $B$ is.
\end{exa}

\begin{cor} \label{a7.8}
Let $A$ be a connected graded $k$-algebra with a balanced
dualizing complex. Suppose that every graded prime quotient
$A / \mfrak{p}$ satisfies one of the following conditions:
\begin{enumerate}
\rmitem{i} $A / \mfrak{p}$ is graded Auslander \tup{(}resp.\ and
graded $\opn{Kdim}$-\tup{(}or $\opn{GKdim}$\tup{)}-Macaulay\tup{)}; or
\rmitem{ii} $A / \mfrak{p}$ has a normal element of positive degree.
\end{enumerate}
Then $A$ is graded Auslander \tup{(}resp.\ and
graded $\opn{Kdim}$-\tup{(}or $\opn{GKdim}$-\tup{)}Macaulay\tup{)}
\end{cor}

\begin{proof}
Use Proposition \ref{a6.16}, Theorem \ref{a7.1} and
induction on $\opn{Kdim} A$.
\end{proof}

Recall that $A$ has {\em enough normal elements} if every graded prime
quotient $A / \mfrak{p}$ (except for $\mfrak{p} = \mfrak{m}$)
contains a normal element of positive degree. By Corollary \ref{a7.8},
a noetherian connected graded algebra with enough normal elements is
Auslander and $\opn{GKdim}$-Macaulay. In the rest of this section we
generalize this statement.

\begin{cor} \label{a7.9}
Let $A$ and $B$ be noetherian connected graded $k$-algebras.
Assume $A$ has a balanced dualizing complex and is graded Auslander
\tup{(}resp.\ and graded $\opn{Kdim}$\tup{(}or
$\opn{GKdim}$\tup{)}-Macaulay\tup{)}, and $B$ has
enough normal elements. Then $A \otimes B$ is a noetherian connected
graded $k$-algebra with a balanced dualizing complex, and it is
graded Auslander \tup{(}resp.\ and
graded $\opn{Kdim}$\tup{(}or $\opn{GKdim}$\tup{)}-Macaulay\tup{)}.
\end{cor}

\begin{proof}
We first prove that $A\otimes B$ is noetherian by induction on
$\opn{Kdim} B$. Suppose $\{\mfrak{p}_1,\cdots, \mfrak{p}_n\}$ is
the complete set of minimal graded primes of $B$. Then some product
$\prod_{i=1}^t \mfrak{p}_{n_i}$ is zero.
For any $s \geq 1$ let
$W_{s} := \prod_{i = 1}^{s - 1} \mfrak{p}_{n_{i}} /
\prod_{i = 1}^{s} \mfrak{p}_{n_{i}}$,
where $W_{1}= B / \mfrak{p}_{n_{1}}$, and let
$B_{s} = B / \mfrak{p}_{n_{s}}$.
Since $W_{s}$ is a finitely generated $B_{s}$-module,
$A \otimes W_{s}$ is a finitely generated $A \otimes B_{s}$-module.
Hence it suffices to show that each $A \otimes B_{s}$
is noetherian. This reduces to the case when $B$ is prime. Hence
we may assume $B$ has a regular normal element $t$ of positive degree.
It is obvious that $1 \otimes t$ is a regular normal element in
$A \otimes B$. By induction hypothesis, $A \otimes B / (t)$ is
noetherian, and therefore $A \otimes B$ is noetherian by \cite{ATV}
Theorem 8.2.

The graded Auslander and Macaulay properties follow from the same
inductive procedure and Theorem \ref{a7.1}.
\end{proof}

In \cite{SZ} Theorem 3.10 it was shown that a connected graded PI
algebra of finite injective dimension (i.e.\ a Gorenstein algebra)
is graded Auslander-Gorenstein and graded $\opn{GKdim}$-Macaulay.
This was extended in \cite{Zh} to an algebra $A$ having
enough normal elements. Corollary \ref{a7.9} (when $A=k$)
generalizes these theorems by eliminating the Gorenstein condition.
We extend the result further in Theorems \ref{a7.11} and \ref{a7.12}
below.

\begin{dfn} \label{a7.10}
Let $R$ be a balanced dualizing complex over $A$.
\begin{enumerate}
\item We say two graded $A$-modules $M$ and $N$ are {\em similar} if
there are isomorphisms $M \cong N$ and
$\mrm{R} \opn{Hom}^{\mrm{gr}}_{A}(M, R) \cong
\mrm{R} \opn{Hom}^{\mrm{gr}}_{A}(N, R)$
in $\msf{D}(\cat{GrMod} k)$.
\item The algebra $A$ satisfies the {\em similar submodule condition}
if every nonzero, $\mfrak{m}$-torsion-free, finitely generated, graded
$A$-module $M$ has
graded submodules $N' \subsetneqq N \subset M$ with
$N'$ similar to $N(-l)$ for some $l > 0$.
\end{enumerate}
\end{dfn}

We remark that two complexes
$M, N \in \msf{D}^{\mrm{b}}_{\mrm{f}}(\cat{GrMod} A)$
are isomorphic in \linebreak
$\msf{D}(\cat{GrMod} k)$ if and only if they have equal
Hilbert functions, namely
$\opn{rank}_{k} (\mrm{H}^{i} M)_{j} =
\opn{rank}_{k} (\mrm{H}^{i} N)_{j}$
for all $i, j \in \mbb{Z}$.
This definition of similar submodule condition is equivalent to
the definition given in \cite{Zh} Section 2 when $A$ is AS-Gorenstein,
which is the case considered there. Also, as mentioned in \cite{Zh}
Section 2, there are algebras which do not satisfy the similar
submodule condition.

\begin{thm} \label{a7.11}
Assume $A$ is a noetherian connected graded $k$-algebra and satisfies
one of the following:
\begin{enumerate}
\rmitem{i} $A$ is a PI algebra.
\rmitem{ii} $A$ is graded-FBN.
\rmitem{iii} $A$ has enough normal elements.
\end{enumerate}
Then $A$ has a balanced dualizing complex, and the
similar submodule condition holds.
\end{thm}

Note that if $A$ is PI, then $A$ is graded FBN and has
enough normal elements.

\begin{proof}
The existence of a balanced dualizing complex follows
Theorem \ref{a6.5}, together with \cite{AZ} Theorems 8.8 and 8.13.
The statement about the similar submodule condition is proved
like in \cite{Zh} Section 2 (first paragraph and Proposition 2.3).
\end{proof}

We now prove a generalized version of \cite{Zh} Theorem 3.1.

\begin{thm} \label{a7.12}
Let $A$ be a noetherian connected graded $k$-algebra
which has a balanced dualizing complex $R$ and
satisfies the similar submodule condition.
Then
\begin{enumerate}
\item $A$ is graded Auslander.
\item $A$ is graded $\opn{Kdim}$- and $\opn{GKdim}$-Macaulay.
More precisely
\[ \opn{Cdim}_{R} M = \opn{Kdim} M = \opn{GKdim} M \]
for every finitely generated graded left or right $A$-module $M$.
\end{enumerate}
\end{thm}

\begin{proof}
First we observe that \cite{Zh} Lemma 2.2 holds (same proof), and
hence
$\opn{Kdim} M$ \linebreak
$\geq \opn{GKdim} M$ and $\opn{GKdim} M < \infty$
for every finitely generated graded $A$-module $M$. Therefore,
replacing ${\opn{Ext}}^{i}(-, A)$ with
$\opn{Ext}^{i}_{A}(-, R)$, the proof of \cite{Zh} Theorem 3.1 can be
copied verbatim. Let us just mention the key point of the proof.
We prove by induction on $\opn{GKdim} M$ that for every
finitely generated graded $A$-module (resp.\  $A^{\circ}$-module) $M$:
\begin{enumerate}
\rmitem{a} $j_{R} (M) = - \opn{GKdim} M$.
\rmitem{b} $\opn{GKdim} \opn{Ext}^{j}_{A}(M, R) =
\opn{GKdim} M$,
where $j = j_{R}(M)$.
\rmitem{c} For every $j_{R}(M) \leq i \leq 0$ one has
$\opn{GKdim} \opn{Ext}^{i}_{A}(M, R) \leq -i$.
\end{enumerate}
This implies that $A$ is graded Auslander and graded
$\opn{GKdim}$-Macaulay. By \cite{Zh} Lemma 2.2,
$\opn{Kdim} M \geq \opn{GKdim} M= \opn{Cdim} M$, and by Theorem
\ref{a6.11},
$\opn{Kdim} M \leq \opn{Cdim} M$.
Hence $A$ is graded $\opn{Kdim}$-Macaulay.
\end{proof}

Theorem \ref{a0.6} is an immediate consequence of
Theorems \ref{a7.11} and \ref{a7.12}.
Another consequence is Proposition \ref{a0.2} which we now prove.

\begin{proof}[Proof of Proposition \tup{\ref{a0.2}}]
If $A$ is connected graded, i.e., $A_0=k$, the assertion follows
from Theorems \ref{a7.11} and \ref{a7.12}.
Now assume $A$ is not connected. Let $B := k + A_{\geq 1}$, which is
connected. Since $A$ and $B$
differ by a finite rank $k$-module, $B$ is noetherian.
It remains to verify the following statements.
\begin{enumerate}
\rmitem{a} If $M$ is a finitely generated graded $A$-module, then
$\opn{Kdim}_{A} M = \opn{Kdim}_{B} M$ and
$\opn{GKdim}_{A} M = \opn{GKdim}_{B} M$.
\rmitem{b} If $A$ has enough normal elements, then so does $B$.
\rmitem{c} If $A$ is graded FBN, then so is $B$.
\end{enumerate}

\medskip \noindent (a)
The statement about $\opn{GKdim}$ is obvious because
$\opn{GKdim}$ is determined by the degree of the Hilbert function of
$M$ (see \cite{Zh} Lemma 2.2(1)).
Next we consider $\opn{Kdim}$. Clearly
$\opn{Kdim}_{A} M = 0$ if and only if $M$ is a finitely generated
$k$-module, if and only if $\opn{Kdim}_{B} M = 0$.
For higher dimension we consider the quotient category
$\cat{QGr} A := \cat{GrMod} A / \cat{M}_{0}$,
where $\cat{M}_{0} = \cat{M}_{0}(\opn{Kdim})$
is the localizing subcategory consisting of $\opn{Kdim} = 0$ modules
(torsion modules in the terminology of \cite{AZ}).
For any $M \in \cat{GrMod}_{\mrm{f}} A$ one has
\[ \opn{Kdim}_{A} M =
\opn{Kdim}_{\cat{GrMod} A} M = \opn{Kdim}_{\cat{QGr} A} M + 1 . \]
Now since $A$ and $B$ differ by a finite rank $k$-module,
$\cat{QGr} A$ is equivalent to $\cat{QGr} B$ by \cite{AZ} Proposition
2.5, so
$\opn{Kdim}_{\cat{QGr} A} M = \opn{Kdim}_{\cat{QGr} B} M$.

\medskip \noindent (b)
Let $\mfrak{p}$ be a graded prime ideal of $B$ which is not
$\mfrak{m}_{B} = B_{\geq 1}$.
Then $A \mfrak{p} A$ and $\mfrak{p}$ differ by a finitely generated
$k$-module, and hence
$\opn{GKdim} A / A \mfrak{p} A = \opn{GKdim} B / \mfrak{p}$.
Let $\mfrak{q}$ be a graded prime of $A$ minimal over $A \mfrak{p} A$
such that $\opn{GKdim} A / \mfrak{q} = \opn{GKdim} B / \mfrak{p}$.
Then the map $B / \mfrak{p} \ar A / \mfrak{q}$ is injective, and
bijective in positive degrees. Thus normal elements of
positive degree in $A / \mfrak{q}$ are also normal elements in
$B / \mfrak{p}$.

\medskip \noindent (c)
The proof is similar to (b) and we leave it to the reader.
\end{proof}

\section{Noetherian Connected Filtrations}

In this section we use filtrations to transfer results of Sections
4-5 on connected graded algebras to non-graded algebras.
Throughout the section $A$ denotes a noetherian $k$-algebra.

Suppose a $k$-module $M$ is given an increasing filtration
$F = \{ F_{n} M \}_{n \in \mbb{Z}}$ with $\bigcap_{n} F_{n} M = 0$
and $\bigcup_{n} F_{n} M = M$.
The {\em Rees module} is the graded $k[t]$-module
\[ \opn{Rees}^{F} M := \bigoplus_{n \in \mbb{Z}} F_{n} M \cdot t^{n}
\subset M[t, t^{-1}] \]
where $t$ is a central indeterminate. It is easy to check that
\[ (\opn{Rees}^{F} M) / (t - 1) \cdot (\opn{Rees}^{F} M) \cong M \]
and
\[ (\opn{Rees}^{F} M) / t \cdot (\opn{Rees}^{F} M) \cong
\opn{gr}^{F} M = \bigoplus_{n \in \mbb{Z}} F_{n} M / F_{n - 1} M . \]

Now let
$F = \{ F_{n} A \}_{n \in \mbb{Z}}$ be a filtration of $A$ such that
$F_{n} A \cdot F_{m} A \subset F_{n + m} A$.
The graded $k$-algebra $\opn{Rees}^{F} A$ is called the
{\em Rees algebra}.

\begin{dfn}
If the Rees algebra $\opn{Rees}^{F} A$ is a noetherian connected
graded $k$-algebra then $F$ is called a {\em noetherian
connected filtration}.
\end{dfn}

Observe that $\opn{Rees}^{F} A$ is connected graded iff
$F_{0} A \cong k$, $F_{-1} A = 0$ and
$\opn{rank}_{k} F_{n} A$ \linebreak $< \infty$.
If $\opn{Rees}^{F} A$ is noetherian then so is the associated graded
algebra $\opn{gr}^{F} A$. By \cite{ATV} Theorem 8.2 the converse is
also true -- if $\opn{gr}^{F} A$ is noetherian then so is
$\opn{Rees}^{F} A$.

A filtration $\{ F_{n} A \}$ on an $A$-module $M$ with
$F_{n} A \cdot F_{m} M \subset F_{n + m} M$
gives a graded module
$\opn{Rees}^{F} M$ over $\opn{Rees}^{F} A$.
We say $\{ F_{n} A \}$ is a {\em good filtration} if
$\opn{Rees}^F M$ is a finitely generated $(\opn{Rees}^F A)$-module.

The main result of this section is the following theorem. Part 1
is due to Van den Bergh (\cite{VdB} Theorem 8.6), but there was a
subtle flaw in his statement: the shift by $-1$ was missing.

\begin{thm} \label{a8.3}
Let $F$ be a noetherian connected filtration on $A$ and let
$\tilde{A} := \opn{Rees}^{F} A$.
\begin{enumerate}
\item If $\tilde{A}$ has a balanced dualizing complex $\tilde{R}$ then
\[ R := (A \otimes_{\tilde{A}} \tilde{R} \otimes_{\tilde{A}} A)
[-1] \]
is a rigid dualizing complex over $A$.
\item If $\tilde{R}$ is graded Auslander then $R$ is Auslander.
\item Suppose $\tilde{R}$ is graded Auslander.
If $\tilde{R}$ is graded $\opn{GKdim}$-Macaulay then $R$ is
$\opn{GKdim}$-Macaulay.
\end{enumerate}
\end{thm}

The proof comes after this lemma.

Consider the functor $\pi : \cat{GrMod} k[t] \ar \cat{Mod} k$,
$\tilde{M} \mapsto \tilde{M} / (t - 1) \tilde{M}$.
Write $\tilde{A}_{t} := A[t,t^{-1}]$.

\begin{lem} \label{a8.1}
\begin{enumerate}
\item The functor $\pi$ is exact.
\item If $\tilde{I} \in \cat{GrMod} \tilde{A}$ is injective then
$\pi \tilde{I} \in \cat{Mod} A$ is injective.
\item There is a functorial isomorphism
\[ \mrm{R} \opn{Hom}_{A}(\pi \tilde{M}, \pi \tilde{N}) \cong
\pi \mrm{R} \opn{Hom}^{\mrm{gr}}_{\tilde{A}}(\tilde{M}, \tilde{N}) \]
for
$\tilde{M}, \tilde{N} \in \msf{D}^{\mrm{b}}_{\mrm{f}}(\cat{GrMod}
\tilde{A})$.
\end{enumerate}
\end{lem}

\begin{proof}
1, 2. $\pi$ is a composition of the functors ``localization by $t$''
$(-)_{t} : \cat{GrMod} \tilde{A} \ar \cat{GrMod} \tilde{A}_{t}$
and ``taking degree $0$ component''
$(-)_{0} : \cat{GrMod} \tilde{A}_{t} \ar \cat{Mod} A$,
both of which are exact. The second is even an equivalence, so
injectives go to injectives. Since $\tilde{A}_{t}$ is noetherian a
standard argument shows that the $\tilde{A}_{t}$-module
$\tilde{I}_{t}$ is injective.

\medskip \noindent 2.\
There is a functorial morphism
$\psi :
\pi \mrm{R} \opn{Hom}^{\mrm{gr}}_{\tilde{A}}(\tilde{M}, \tilde{N})
\ar \mrm{R} \opn{Hom}_{A}(\pi \tilde{M}, \pi \tilde{N})$.
Fixing $\tilde{N}$ these are way-out right contravariant functors of
$\tilde{M}$. By \cite{RD} Proposition I.7.1(iv) -- reversed --
it's enough to check that $\psi$ is an isomorphism when
$\tilde{M} = \tilde{A}$, and then it's trivial.
\end{proof}

\begin{proof}[Proof of Theorem \tup{\ref{a8.3}}]
1.\ By Proposition \ref{a6.25} and the graded version of Corollary
\ref{a4.5} each cohomology module
$\mrm{H}^{q} \tilde{R}$ is $k[t]$-central (cf.\ Remark \ref{a8.20}
below). Define
\[ \tilde{R}_{t} := \tilde{A}_{t} \otimes_{\tilde{A}} \tilde{R}
\otimes_{\tilde{A}} \tilde{A}_{t}
\in \msf{D}(\cat{GrMod} (\tilde{A}_{t})^{\mrm{e}}) . \]
We see that the homomorphisms
$\tilde{A}_{t} \otimes_{\tilde{A}} \tilde{R} \ar \tilde{R}_{t}$
and
$\tilde{R} \otimes_{\tilde{A}} \tilde{A}_{t} \ar \tilde{R}_{t}$
are quasi-isomorphisms.

Because $\tilde{A} \ar A = \pi \tilde{A}$ is flat we also have
\[ A \otimes_{\tilde{A}} \tilde{R}[-1] \cong
\tilde{R}[-1] \otimes_{\tilde{A}} A \cong R
\in \msf{D}(\cat{Mod} \tilde{A}^{\mrm{e}}) . \]
Considering only $\tilde{A}$-modules we have
$R \cong \pi \tilde{R}[-1]$, so by Lemma \ref{a8.3}, $R$ has finite
injective dimension and finitely generated cohomologies
over $A$. By symmetry this is true also over $A^{\circ}$.
Part 3 of the lemma implies that
$\mrm{R} \opn{Hom}_{A}(R, R) \cong \pi \tilde{A} = A$,
and likewise
$\mrm{R} \opn{Hom}_{A^{\circ}}(R, R) = A$.
We conclude that $R$ is dualizing.

Now let's prove $R$ is rigid.
There is a functorial isomorphism
$\pi \tilde{M} \cong (\tilde{M}_{t})_{0}$ for
$\tilde{M} \in \cat{GrMod} k[t]$, and the algebra $A^{\mrm{e}}$ is
$\mbb{Z}^{2}$-graded. Therefore we get
\[ A^{\mrm{e}} \otimes_{\tilde{A}^{\mrm{e}}} (
\tilde{R} \otimes \tilde{R}) \cong
(A \otimes_{\tilde{A}} \tilde{R}) \otimes
(\tilde{R} \otimes_{\tilde{A}} A) \cong
(\tilde{R}_{t})_{0} \otimes (\tilde{R}_{t})_{0} \cong
(\tilde{R}_{t} \otimes \tilde{R}_{t})_{(0, 0)} . \]
The algebra $(\tilde{A}_{t})^{\mrm{e}}$ is strongly $\mbb{Z}$-graded,
and its degree $0$ component is
\[ (\tilde{A}_{t})^{\mrm{e}}_{0} \cong A^{\mrm{e}}[s, s^{-1}] \cong
k[s, s^{-1}] \otimes A^{\mrm{e}} \]
where
$s := t \otimes t^{-1}$.
Applying
$(\tilde{A}_{t})^{\mrm{e}} \otimes_{\tilde{A}^{\mrm{e}}} -$
to
\[ \tilde{R} \cong
\mrm{R} \opn{Hom}^{\mrm{gr}}_{\tilde{A}^{\mrm{e}}}(\tilde{A},
\tilde{R} \otimes \tilde{R}) \]
we obtain
\[ \tilde{R}_{t} \cong
\mrm{R} \opn{Hom}^{\mrm{gr}}_{(\tilde{A}_{t})^{\mrm{e}}}
(\tilde{A}_{t}, \tilde{R}_{t} \otimes \tilde{R}_{t}) \]
and taking degree $0$ components we get
\[ \begin{aligned}
R[1] & \cong
\mrm{R} \opn{Hom}^{\mrm{gr}}_{(\tilde{A}_{t})^{\mrm{e}}_{0}}
(A, (\tilde{R}_{t} \otimes \tilde{R}_{t})_{0}) \\
& \cong
\mrm{R} \opn{Hom}_{k[s, s^{-1}] \otimes A^{\mrm{e}}}
(A, k[s, s^{-1}] \otimes (R[1] \otimes R[1])) .
\end{aligned} \]
But
\[ k \cong \left(k[s, s^{-1}] \xrightarrow{s - 1} k[s, s^{-1}]\right)
\in \msf{D}(\cat{Mod} k[s, s^{-1}]) , \]
and $A$ is a finitely presented $A^{\mrm{e}}$-module, so
\[ \mrm{R} \opn{Hom}_{k[s, s^{-1}] \otimes A^{\mrm{e}}}
(k \otimes A, k[s, s^{-1}] \otimes M) \cong
\mrm{R} \opn{Hom}_{A^{\mrm{e}}} (A, M) [-1] \]
for any complex of $A^{\mrm{e}}$-modules $M$.

\medskip \noindent 2.\
Let $\tilde{M}$ be a $t$-torsion-free finitely generated
graded $\tilde{A}$-module. We claim that
\begin{equation} \label{eqn7.4}
j_{R; A} (\pi \tilde{M}) = j_{\tilde{R}; \tilde{A}} \tilde{M} + 1 .
\end{equation}
First by Lemma \ref{a8.1} we get for any $q$:
\begin{equation} 
\opn{Ext}^{-q + 1}_{A}(\pi \tilde{M}, R) \cong
\pi \opn{Ext}^{-q}_{\tilde{A}}(\tilde{M}, \tilde{R}) ,
\end{equation}
so
$j_{R; A} (\pi \tilde{M}) \geq j_{\tilde{R}; \tilde{A}} \tilde{M} + 1$.
Now take
$q := -j_{\tilde{R}; \tilde{A}} \tilde{M}$
and write
$\tilde{N} := \opn{Ext}^{-q}_{\tilde{A}}(\tilde{M}, \tilde{R})$.
If
$\pi \tilde{N} = 0$ then $t^{l} \tilde{N} = 0$ for some $l > 0$.
But in the Ext spectral sequence converging to $\tilde{M}$
(see proof of Theorem \ref{a2.19}) the dominant term is
$\opn{Ext}^{-q}_{\tilde{A}^{\circ}}(\tilde{N}, \tilde{R})$
which is killed by $t^{l}$. We get
$\opn{Cdim}_{\tilde{R}} t^{l} \tilde{M} <
\opn{Cdim}_{\tilde{R}} \tilde{M} = q$,
which is absurd since
$t^{l} \tilde{M} \cong \tilde{M}(-l)$.

Given finitely generated $A$-modules $M \subset N$, take any
good filtration $\{ F_{n} M \}$, and let
$F_{n} N := N \cap F_{n} M$, $\tilde{M} := \opn{Rees}^{F} M$ and
$\tilde{N} := \opn{Rees}^{F} N$.
Since $\tilde{M}$ and $\tilde{N}$ are $t$-torsion-free we see that
\[ \opn{Cdim}_{R} M = \opn{Cdim}_{\tilde{R}} \tilde{M} - 1
\leq \opn{Cdim}_{\tilde{R}} \tilde{N} - 1 = \opn{Cdim}_{R} N . \]

Finally let $L := \opn{Ext}^{-q}_{A}(M, R)$. Then, with
$\tilde{M} := \opn{Rees}^{F} M$ and
$\tilde{L} :=$ \linebreak
$\opn{Ext}^{q - 1}_{\tilde{A}}(\tilde{M}, \tilde{R})$,
we see that
$j_{R; A} L \geq j_{\tilde{R}; \tilde{A}} \tilde{L} + 1 \geq -q$,
verifying the Auslander condition on one side. By symmetry it holds
also on the other side.

\medskip \noindent 3.\
By the proof of part 2, given a finitely generated $A$-module $M$ one
has
$\opn{Cdim}_{R; A} M$ \linebreak
$+ 1 = \opn{Cdim}_{\tilde{R}; \tilde{A}} \tilde{M}$
where
$\tilde{M} := \opn{Rees}^{F} M$ w.r.t.\ any good filtration
$\{ F_{n} M \}$.
But because
$\tilde{M}_{t} \cong M \otimes k[t, t^{-1}]$
we also have
\[ \opn{GKdim}_{A} M + 1 =
\opn{GKdim}_{\tilde{A}_{t}} \tilde{M}_{t} =
\opn{GKdim}_{\tilde{A}} \tilde{M} . \]
\end{proof}

\begin{rem} \label{a8.20}
It is not too hard to show that $R$ could be chosen to be a
$k[t]$-central complex of graded $A$-bimodules.
\end{rem}

Observe that the rigid dualizing complex $R$ has
$\mrm{H}^{q} R = 0$ for $q > 0$, since $\mrm{H}^{0} \tilde{R}$
is $t$-torsion and $\mrm{H}^{q} \tilde{R} = 0$ for $q > 0$.

\begin{cor} \label{a8.6}
Suppose $A$ is connected graded and $R$ is a balanced dualizing
complex over $A$. Then $R$ is a rigid dualizing complex over $A$
in the ungraded sense.
\end{cor}

\begin{proof}
First let us note that for any graded $A$-module $M$ we can define a
filtration
$F_{n} M := \bigoplus_{i \leq n} M_{i}$.
Then we have a functorial isomorphism of graded $A[t]$-modules
$\opn{Rees}^{F} M \iso M[t]$.

In particular we get $\tilde{A} \cong A[t]$, so by Example \ref{a7.7}
we know that
$\tilde{R} = R[t](-1)[1]$
is the balanced dualizing complex of $\tilde{A}$. But then
$\pi \tilde{R}[-1] = \pi R[t](-1) \cong R$
in $\msf{D}(\cat{Mod} A^{\mrm{e}})$, and this is rigid by
Theorem \ref{a8.3}(1).
\end{proof}

The next corollary implies Theorem \ref{a0.5}.

\begin{cor} \label{a8.8}
Suppose $A$ has a noetherian connected filtration $F$,
and let $\bar{A} := \opn{gr}^F A$.
\begin{enumerate}
\item If $\bar{A}$ has a balanced dualizing complex $\bar{R}$, then
$A$ has a rigid dualizing complex $R$.
\item If $\bar{R}$ is graded Auslander then $R$ is Auslander.
\item If $\bar{R}$ is also graded $\opn{GKdim}$-Macaulay, then $R$ is
$\opn{GKdim}$-Macaulay.
\end{enumerate}
\end{cor}

\begin{proof}
According to Theorem \ref{a7.1}, the Rees algebra
$\tilde{A} = \opn{Rees}^{F} A$ inherits these properties from $\bar{A}$.
And by Theorem \ref{a8.3} they pass to $A$.
\end{proof}

\begin{cor} \label{a8.9}
Suppose $A$ has a noetherian connected filtration $F$,
and $\bar{A} := \opn{gr}^{F} A$ satisfies either of the following:
\begin{enumerate}
\rmitem{i} $\bar{A}$ is a PI algebra.
\rmitem{ii} $\bar{A}$ is graded-FBN.
\rmitem{iii} $\bar{A}$ has enough normal elements.
\end{enumerate}
Then $A$ has an Auslander, $\opn{GKdim}$-Macaulay, rigid dualizing
complex.
\end{cor}

\begin{proof}
Combine Corollary \ref{a8.8} and
Theorems \ref{a7.11} and \ref{a7.12}.
\end{proof}

Here are some examples of algebras which admit noetherian connected
filtrations.

\begin{exa} \label{a8.10}
If $A$ is a noetherian connected graded algebra, then the filtration
$F_{n} A = \bigoplus_{i \leq n} A_{i}$ is a noetherian connected
filtration.
\end{exa}

\begin{exa} \label{a8.26}
Suppose $A$ is generated by elements $x_{1}, \ldots, x_{n}$, and for
every $i \neq j$ there is some relation
\begin{equation} \label{a8.27}
x_{j} x_{i} = q_{i, j} x_{i} x_{j} + a_{i, j} x_{i} + b_{i, j} x_{j}
+ c_{i, j}
\end{equation}
with
$q_{i, j}, a_{i, j}, b_{i, j}, c_{i, j} \in k$.
Let
$V := k + \sum k \cdot x_{i} \subset A$ and define a filtration
$F_{n} A := V^{n}$. Then $\opn{gr}^{F} A$ is a quotient of the skew
polynomial algebra $k_{q}[x_{1}, \ldots, x_{n}]$, so $F$ is a
noetherian connected filtration. Furthermore $\opn{gr}^{F} A$ has
enough normal elements (namely the $x_{i}$); so Corollary \ref{a8.9}
holds. It is easy to check that the relations (\ref{a8.26}) are
satisfied in the following classes of algebras:
\begin{enumerate}
\rmitem{i} Commutative affine algebras.
\rmitem{ii} Weyl algebras, enveloping algebras of finite
dimensional Lie algebras and their quotients.
\rmitem{iii} Most classes of quantum algebras listed in \cite{GL}.
\end{enumerate}
\end{exa}

Recall that a homomorphism $f : A \ar B$ is called finite if $B$ is a
finite left and right $A$-module. $f$ is {\em centralizing} if
$B = A \cdot \mrm{Z}_{B}(A)$. Thus $f$ is finite centralizing iff
there exist $b_{1}, \ldots, b_{m} \in \mrm{Z}_{B}(A)$ such that
$B = \sum A \cdot b_{i}$.

\begin{lem} \label{a8.13}
Suppose $f : A \ar B$ is a finite centralizing homomorphism and $F$ is
a noetherian connected filtration on $A$. Then there is a noetherian
connected filtration $F$ on $B$ such that $f$ preserves the
filtrations and
$\opn{Rees}^{F}(f) : \opn{Rees}^{F} A \ar \opn{Rees}^{F} B$
is finite.
\end{lem}

\begin{proof}
Let $b_{1}, \ldots, b_{m}$ be elements of $B$ which commute
with $A$ and $B = \sum A \cdot b_{i}$. Choose elements
$a_{i,j,l} \in A$ such that  $b_{i}b_j = \sum_{l} a_{i,j,l} b_{l}$.
Let $n_{0} > 0$ be large enough such that all $a_{i,j,l}$ are in
$F_{n_{0}} A$. Define
\[  F_{n} B := F_{n} A \cdot 1 +
\sum_{i} F_{n - n_{0}} A \cdot b_{i} \subset B . \]
Clearly this is a connected filtration. Since the elements
$(1, b_{1}, \ldots, b_{m})$ determine a surjective bimodule
homomorphism
\[ \opn{Rees}^{F} A \oplus (\opn{Rees}^{F} A)^{m}(-n_{0})
\surj \opn{Rees}^{F} B \]
we see that $\opn{Rees}^{F} B$ is noetherian.
\end{proof}

\begin{exa}
If $A$ is an affine $k$-algebra finite over its center, then there is
a finite centralizing homomorphism
$k[t_{1}, \ldots, t_{n}] \ar A$ from a commutative polynomial
algebra. By the lemma and Example \ref{a8.10}, $A$ has a
noetherian connected filtration.
Thus $A$ satisfies the assumptions of Corollary \ref{a8.9}.
\end{exa}

\begin{exa} \label{exa8.1}
Here is an example of a prime PI algebra $A$ which is not finite
over its center yet has a noetherian connected filtration
(Schelter's Example, \cite{Ro} p.\ 492 Exercise 27).
Let $t, \lambda_{1}, \lambda_{2}$ be commuting
indeterminates of degree $1$. Define
\[ \begin{aligned}
\tilde{C} & := \mbb{Q}[\sqrt{2}, \sqrt{3}, \lambda_{1}, \lambda_{2},
t]\\
\tilde{C}_{1} & := \mbb{Q}[\sqrt{6}, t \sqrt{2} + \lambda_{1},
\lambda_{2}, \lambda_{2} \sqrt{2}, t]\\
\tilde{C}_{2} & := \mbb{Q}[\sqrt{6}, t \sqrt{3} + \lambda_{1},
\lambda_{2}, \lambda_{2} \sqrt{2}, t]\\
\tilde{M} & := \tilde{C}_{1} \lambda_{2} +
\tilde{C}_{1} \lambda_{2} \sqrt{2} \subset \tilde{C} \\
\tilde{A} & := \left[ \begin{matrix}
\tilde{C}_{1} & \tilde{M} \\
\tilde{M} & \tilde{C}_{2} \end{matrix} \right] .
\end{aligned} \]
Then $\tilde{A}$ is a noetherian graded algebra and $\tilde{A}_{0}$ is
finite over $k := \mbb{Q}$. The quotient
$A := \tilde{A} / (t - 1)$ acquires a filtration
$\{ F_{n} A \}$, and if we modify it by setting $F_{0} A := \mbb{Q}$
this becomes a connected filtration.
But $A$ is not finite over its center.
\end{exa}

\begin{que}
Does every noetherian affine PI $k$-algebra admit a noetherian
connected filtration?
This seems to be a hard question. A similar one was posed by M. Lorenz 
over ten years ago (see \cite{Lo} p.\ 436).
\end{que}

In Section 3 we found that rigid dualizing complexes are
sometimes functorial w.r.t.\ finite algebra homomorphisms, via the
trace morphism. Here is such an instance:

\begin{thm} \label{a8.14}
Let $A \ar B$ be a finite centralizing homomorphism.
Suppose $A$ has a noetherian connected filtration $F$ and
$\opn{gr}^{F} A$ has a balanced dualizing complex. Then
$A$ and $B$ have rigid dualizing complexes $R_{A}$ and $R_{B}$
respectively, and the trace morphism
$\opn{Tr}_{B / A} : R_{B} \ar R_{A}$ 
of Definition \tup{\ref{a4.7}} exists.
\end{thm}

\begin{proof}
By Lemma \ref{a8.13} we get a finite homomorphism of graded algebras
$\tilde{A} = \opn{Rees}^{F} A \ar \tilde{B} = \opn{Rees}^{F} B$.
So according to Theorem \ref{a6.14} the trace morphism
$\opn{Tr}_{\tilde{B} / \tilde{B}} : R_{\tilde{B}} \ar R_{\tilde{A}}$
exists. Now apply the functor $\pi$ and use Theorem \ref{a8.3}.
\end{proof}

For applications of this result see Proposition \ref{a4.9}.

Let $\sigma$ be a $k$-algebra automorphism $A$. Recall that
$A^{\sigma}$ is the invertible bimodule with basis $e$
satisfying $e \cdot a = \sigma(a) \cdot e$. A noetherian
connected graded $k$-algebra $B$ is called {\em AS-Gorenstein} if $B$
satisfies $\chi$ and the bimodule $B$ has finite injective
dimension on both sides. Here ``AS'' stands for Artin-Schelter.
We say $B$ is {\em AS-regular} if $B$ is AS-Gorenstein
and $\opn{gl.dim} B <\infty$.

\begin{prop} \label{a8.15}
Suppose $A$ has a noetherian connected filtration $F$
and $\opn{gr}^{F} A$ is AS-Gorenstein \tup{(}resp. AS-regular\tup{)}.
Then the following statements hold.
\begin{enumerate}
\item $A$ is a Gorenstein \tup{(}resp. regular\tup{)} algebra.
\item The rigid dualizing complex of $A$ is
$R = A^{\sigma}[n]$ where $n$ is an integer and $\sigma$
is some $k$-algebra automorphism of $A$.
\item If $\opn{gr}^{F} A$ is graded Auslander \tup{(}resp.\ and
graded $\opn{GKdim}$-Macaulay\tup{)}, then $A$ is Auslander-Gorenstein
\tup{(}resp. and Cohen-Macaulay\tup{)} in the sense of \cite{Bj}.
\item Let $B = A / \mfrak{a}$ be any quotient algebra, $M$
any $B$-module and $q$ an integer. Then the twisted module
$\opn{Ext}^{q}_{A}(M, A)^{\sigma}$ is a
$B^{\circ}$-module.
\end{enumerate}
\end{prop}

\begin{proof}
1, 2. Let $n$ be the injective dimension of $\opn{gr}^{F} A$, and let
$\tilde{A} := \opn{Rees}^{F} A$. By \cite{Lev} Theorem 3.6
the injective dimension of $\tilde{A}$ is $n + 1$.
Since $\tilde{A}$ satisfies $\chi$ it is AS-Gorenstein.
So the balanced dualizing complex of $\tilde{A}$ is
$\tilde{R} = \tilde{A}^{\tilde{\sigma}}(d)[n + 1]$ for some graded
automorphism $\tilde{\sigma}$ and for some integer $d$.
By Theorem \ref{a8.3} the rigid dualizing complex of $A$ is
$R = \pi \tilde{A}^{\tilde{\sigma}}[n]$.
Since $\tilde{A}^{\tilde{\sigma}} = \mrm{H}^{-n - 1} \tilde{R}$
this is $k[t]$-central; so $\tilde{\sigma}(t) = t$.
We see there is an induced automorphism $\sigma$ of $A$ and
$R = A^{\sigma}[n]$.

If $\opn{gr} A$ has finite global dimension, then so does $A$.

\noindent 3.
This follows from 2 and Theorems \ref{a7.1} and \ref{a8.3}.

\noindent 4.
It follows from Theorem \ref{a8.14} and Proposition \ref{a4.9}(1) 
that
\[ \opn{Ext}^{q}_{A}(M, A)^{\sigma} \cong
\opn{Ext}^{q - n}_{A}(M, A^{\sigma}[n]) \cong
\opn{Ext}^{q - n}_{A}(M, R) \cong
\opn{Ext}^{q - n}_{B}(M, R_{B}) \]
where $R_{B}$ is the rigid dualizing complex of $B$.
\end{proof}

Example \ref{a2.2} shows that Proposition \ref{a8.15}(2,3) might fail
even if $A$ is Gorenstein. The next example shows that $\sigma$
could be nontrivial.

\begin{exa} \label{a8.16}
Let $A$ be the quantum plane
$k_q[x,y] := k \langle x, y \rangle / (y x - q x y)$
for $q \in k^{\times}$ with $q^{2} \neq 1$.
The automorphism $\sigma$ in Proposition
\ref{a8.15}(4) is
$\sigma(x) = qx$, $\sigma(y) = q^{-1}y$
(cf.\ \cite{Ye1} Examples 6.21 and 7.14).
Consider the ideal $I = A \cdot f \cdot A$ where $f := x - y$,
which is not normal. An easy computation shows that
$B \cong k[\epsilon]$ with $\epsilon^{2} = 0$, and
$\epsilon \equiv x \equiv y \ (\opn{mod} I)$.
Now consider the graded
$A^{\mrm{e}}$-module $N := \opn{Ext}^{2}_{A}(B, A)$.
One has
$N^{\sigma}(-2) \cong B^{*}$
as $A^{\mrm{e}}$-modules, so $N$ is killed by
$q x - q^{-1} y = \sigma(f) \in A^{\circ}$,
and hence $N$ can not be a $B^{\circ}$-module.
\end{exa}

\begin{exa} \label{a8.17}
Let $A$ be the Weyl algebra
$k \langle x, y \rangle / (x y - y x - 1)$.
Take the standard
filtration $F_n A=(k+kx+ky)^n$. Then the Rees algebra $\tilde{A}$
is generated by $x, y, t$ with $t$ central and $x y - y x = t^{2}$.
$\tilde{A}$ is an Artin-Schelter regular algebra of global dimension
$3$, so its balanced dualizing complex is
$\tilde{A}^{\tilde{\sigma}}(-3)[3]$ for some automorphism
$\tilde{\sigma}$.
Let
$\tilde{B} := \tilde{A} / (t^{2})$, which is a commutative
AS-Gorenstein algebra. As in \cite{Ye1} Theorem 7.18 we find that
$\tilde{\sigma} = 1$. Therefore the rigid dualizing complex of $A$ is
$A[2]$. Observe that $\opn{Cdim} A = 2 = \opn{GKdim} A$.
In \cite{Ye4} we prove the more general statement that if $C$ is any
smooth integral commutative $k$-algebra of dimension $n$, 
$\opn{char} k = 0$ and $A := \mcal{D}(C)$ is the ring
of differential operators, then the rigid dualizing complex of $A$ is
$A[2n]$.
\end{exa}

Suppose $A$ has a noetherian connected filtration $F$. A
{\em two-sided good filtration} on a bimodule $M$ is a filtration
$\{ F_{n} M \}$ such that
$F_{n} A \cdot F_{m} M \subset F_{n + m} M$,
$F_{n} M \cdot F_{m} A \subset F_{n + m} M$
and $\opn{Rees}^{F} M$ is a finitely generated
$(\opn{Rees}^{F} A)$-module on both sides.

\begin{prop} \label{a8.22}
Assume $A$ has a noetherian connected filtration and
$\opn{gr} A$ has a balanced dualizing complex. Let $R$ be the
rigid dualizing complex of $A$. If a bimodule $M$ has a
two-sided good filtration then
\[ \mrm{R} \opn{Hom}_{A}(M, R) \cong
\mrm{R} \opn{Hom}_{A^{\circ}}(M, R) . \]
\end{prop}

\begin{proof}
Let $\tilde{M} := \opn{Rees}^{F} M$. According to Corollary
\ref{a6.15} there is an isomorphism
\[ \mrm{R} \opn{Hom}_{\tilde{A}}(\tilde{M}, \tilde{R}) \cong
\mrm{R} \opn{Hom}_{\tilde{A}^{\circ}}(\tilde{M}, \tilde{R}) \]
in $\cat{D}(\cat{GrMod} A^{\mrm{e}})$,
where $\tilde{R}$ is the balanced dualizing complex of $\tilde{A}$.
Since $\tilde{M}$ is $k[t]$-central we can apply the functor $\pi$.
\end{proof}

Recall the notion of weakly symmetric dimension function
(Definition \ref{a2.22}).

\begin{cor} \label{a8.23}
Assume $A$ has an Auslander rigid dualizing complex, and a
noetherian connected filtration such that
$\opn{gr} A$ has a balanced dualizing complex.
Then $\opn{Cdim}$ is weakly symmetric.
\end{cor}

\begin{proof}
As can be readily verified, the class of $A$-bimodules which
admit two-sided good filtrations is closed under submodules,
quotients and finite direct sums. Given a bimodule $M$ which is a
subquotient of $A$, Proposition \ref{a8.22} applies and hence
$\opn{Cdim}_{A} M = \opn{Cdim}_{A^{\circ}} M$.
\end{proof}

We can now give the

\begin{proof}[Proof of Theorem \tup{\ref{a0.1}}]
By Corollary \ref{a8.9}, $A$ has an Auslander rigid dualizing complex
$R$, and by Corollary \ref{a8.23}, $\opn{Cdim}_{R}$ is weakly
symmetric (this also follows from the $\opn{GKdim}$-Macaulay
property). Now use Theorem \ref{a3.4}.
\end{proof}

The next theorem has the same conclusions as \cite{ASZ} Theorem 6.1,
but our assumptions are much more focused.
Let $\mfrak{n}$ be the prime radical of $A$. Recall that $\mfrak{n}$
is said to be {\em weakly invariant} w.r.t.\ an exact dimension
function $\opn{dim}$ if
$\opn{dim} \mfrak{n} \otimes_{A} M < \opn{dim} A / \mfrak{n} =
\opn{dim} A$
for every finitely generated $A$-module $M$ with
$\opn{dim} M < \opn{dim} A$ (and the same for right modules);
cf.\ \cite{MR} 6.8.13.
A ring is called {\em quasi-Frobenius} if it is artinian and self
injective.

\begin{thm} \label{a8.24}
Let $A$ be an Auslander-Gorenstein noetherian $k$-algebra of injective
dimension $n$. Assume $A$ has a filtration such that
$\opn{gr} A$ is an AS-Gorenstein noetherian connected graded
$k$-algebra. Then
\begin{enumerate}
\item The prime radical $\mfrak{n}$ is weakly invariant.
\item If $\mfrak{p}$ is a minimal prime then
$\opn{Cdim} A / \mfrak{p} = n$.
\item $A$ has a quasi-Frobenius ring of fractions.
\end{enumerate}
\end{thm}

\begin{proof}
By Proposition \ref{a8.15} the rigid dualizing complex of $A$ is
$R = A^{\sigma}[n]$.
According to Corollary \ref{a8.23},
$\opn{Cdim}_{R}$ is weakly symmetric.
Now the function denoted $\delta$ in \cite{ASZ} Theorem 6.1
coincides with $\opn{Cdim}_{R}$, so all assumptions of that theorem
hold.
\end{proof}

The next theorem is due to Gabber in the case when $\opn{gr} A$ is
Auslander-Goren\-stein, and an elegant proof was communicated to us by
Van den Bergh. We extend the result by dropping the Gorenstein
condition.

\begin{thm} \label{a8.18}
Let $A$ be a filtered $k$-algebra such that $\opn{gr} A$ is a
noetherian connected graded $k$-algebra with graded Auslander
balanced dualizing complex. Given a $\opn{Cdim}$-pure $A$-module $M$,
there is a good filtration $\{ F_{n} M \}$ on it such that
$\opn{gr}^{F} M$ is $\opn{Cdim}$-pure.
\end{thm}

\begin{proof}
The basic idea is to start with an arbitrary good filtration on $M$
and to modify it to get purity.

Let $\tilde{A}$ be the Rees algebra of $A$, and let $\tilde{R}$ be
its balanced dualizing complex. So $\tilde{R}$ has the graded
Auslander property.

Let $n := \opn{Cdim} M + 1$.
Choose any good filtration $F'$ on $M$ and let
$\tilde{M} := \opn{Rees}^{F'} M$.
Since $\tilde{M}$ is $t$-torsion-free and
$M = \pi \tilde{M}$, by (\ref{eqn7.4}) we have
$\opn{Cdim} \tilde{M} = n$. If
$\tilde{M}' \subset \tilde{M}$ is any nonzero graded submodule then
because $\pi \tilde{M}' \subset M$,
and because $M$ is pure, we see that
$\opn{Cdim} \tilde{M}' = n$. Thus $\tilde{M}$ is pure.

Set
\[ E(\tilde{M}) =
\opn{Ext}^{-n}_{\tilde{A}^{\circ}}(
\opn{Ext}^{-n}_{\tilde{A}}(M, \tilde{R}), \tilde{R}) . \]
As in Theorem \ref{a2.6}(3) there is an exact sequence
\begin{equation} \label{a8.19}
0 \ar \tilde{M} \ar E(\tilde{M}) \ar
\tilde{Q} \ar 0
\end{equation}
where $\opn{Cdim} \tilde{Q} \leq n - 2$.
Consider the module $\tilde{N}$ which is the $t$-saturation of
$\tilde{M}$ in $E(\tilde{M})$, i.e.\
\[ \tilde{N} := \{ x \in E(\tilde{M})
\mid t^{i} x \in \tilde{M} \text{ for some } i \geq 0 \} . \]
Since $\tilde{N} / \tilde{M}$ is $t$-torsion it follows that
$M = \pi \tilde{M} \ar \pi \tilde{N}$ is bijective; but the
filtration $F$ on $\pi \tilde{N}$ may be different from $F'$.
Observe that $\tilde{N}$ is
$t$-torsion-free, so by lifting back the filtration we obtain
$\tilde{N} \cong \opn{Rees}^{F} \pi \tilde{N}$. By (\ref{a8.19})
we get $\opn{Cdim} \tilde{N} / \tilde{M} \leq n - 2$,
which implies that
$E(\tilde{M}) \ar E(\tilde{N})$ is bijective.
Therefore by changing the good filtration on $M$ from $F'$ to $F$,
we can assume that in (\ref{a8.19})
the module $\tilde{Q}$ is $t$-torsion-free.

Having done so we get a short exact sequence
\[ 0 \ar \pi_{0} \tilde{M} \ar
\pi_{0} E(\tilde{M}) \ar  \pi_{0} \tilde{Q} \ar 0 \]
where $\pi_{0}$ denotes the functor $M \mapsto M / t M$.
Hence in order to prove that
$\opn{gr} M = \pi_{0} \tilde{M}$
is pure it suffices to prove that
$\pi_{0} E(\tilde{M})$ is pure.
Now, using the duality functor
$D = \mrm{R} \opn{Hom}^{\mrm{gr}}_{\tilde{A}}(-, \tilde{R})$
we have
$\tilde{L} := \opn{Ext}^{-n}_{\tilde{A}}(\tilde{M}, \tilde{R}) =
\mrm{H}^{-n} D \tilde{M}$,
which is $t$-torsion-free by Lemma \ref{a7.2}. Therefore by the same
lemma,
$\opn{Cdim} \pi_{0} \tilde{L} = n - 1$.
Now
$E(\tilde{M}) = \mrm{H}^{-n} D^{\circ} \tilde{L}$,
and by formula (\ref{a7.3}) when $i=-n$ we get
\[ \pi_{0} \mrm{H}^{-n} D^{\circ} \tilde{L} \subset
(\mrm{H}^{- (n - 1)} D^{\circ} \pi_{0} \tilde{L})(-1) . \]
But by Theorem \ref{a2.6}, the module
$\mrm{H}^{-(n - 1)} D^{\circ} \pi_{0} \tilde{L}$ is pure of
dimension $n - 1$.
\end{proof}

If $\opn{gr} A$ is a commutative affine $k$-algebra and $M$ is a
finitely generated $A$-module, define $I(M) \subset \opn{gr} A$
to be the prime radical of
$\opn{Ann}_{\opn{gr} A} \opn{gr}^{F} M$ for some good filtration $F$
on $M$. By \cite{MR} Proposition 8.6.17, $I(M)$ is independent of the
choice of good filtration. The {\em characteristic variety} of $M$
is defined to be
\[ \opn{Ch}(M) = \opn{Spec} \opn{gr} A / I(M)  \]
(cf.\ \cite{Co} p.\ 98 for the case when $A$ is a Weyl algebra).

\begin{proof}[Proof of Theorem \tup{\ref{a0.4}}]
Recall that a variety is called pure if all its irreducible
components have the same dimension. The support of a finitely
generated module $N$ over a commutative affine $k$-algebra $B$
is pure iff $N$ is $\opn{Kdim}$-pure. But for a commutative algebra
$\opn{Kdim}_{B} = \opn{GKdim}_{B} = \opn{Cdim}_{B}$.
Here we take $B := \opn{gr} A$ and $N := \opn{gr}^{F} M$.
\end{proof}


\end{document}